\documentclass[12pt]{article}

\usepackage{amsmath,amsthm}
\usepackage{amssymb}
\usepackage[mathscr]{eucal}
\usepackage[usenames]{color}
\usepackage[latin5]{inputenc}
\usepackage{url}
\usepackage{longtable}
\usepackage{graphicx}
\usepackage{epstopdf}
\usepackage{subfigure}
\usepackage{graphics}
\usepackage{amssymb}
\usepackage[mathscr]{eucal}
\usepackage[usenames]{color}
\usepackage{subfigure}
\usepackage{lipsum}
\usepackage{mathtools}
\usepackage{cuted}
\usepackage[
bookmarks=true,
backref=true,
colorlinks=true,
linkcolor=blue,
citecolor=red,
urlcolor=blue]{hyperref}
\newcommand{\om}{\omega}

\newtheorem{lemma}{Lemma}

\newtheorem{theorem}{Theorem}

\newtheorem{remark}{Remark}

\numberwithin{equation}{section}

\def\om{\omega}

\begin{document}

\title{\Large Hopf Bifurcation of a Financial Dynamical System with Delay
}

\author{
Yasemin \c Cal{\i}\c s\thanks{e-mail: calis.ysmn@gmail.com} \quad  Ali Demirci\thanks{e-mail: demircial@itu.edu.tr}\quad  Cihangir \"{O}zemir\thanks{Corresponding author, e-mail: ozemir@itu.edu.tr}\\
\small Department of Mathematics, Faculty of Science and Letters,\\
\small Istanbul Technical University, 34469 Istanbul,
Turkey}

\date{February 21, 2021}

\maketitle

\begin{abstract}
The aim of this work is to investigate the qualitative behaviour of a  financial dynamical system  which contains a time delay. We investigate  the dynamic response of this system of which variables are interest rate, investment demand, price index and average profit margin. We perform a stability analysis at the fixed points and show that the system undergoes a Hopf bifurcation. The bifurcation analyses are supported by numerical simulations.

\vspace{0.5cm}
\end{abstract}

\section{Introduction}
\label{intro}

The interplay between the dynamical systems theory and economic and financial science has been a major subject of research both for the mathematicians and experts of economic fields in the past decades and to date. In dynamical systems literature there are lots of mathematical models related to finance theory.  We would like to mention first the financial dynamical system
		\begin{subequations}\label{eq60}
			\begin{eqnarray}
			&&\dot{x}=z+(y-a)x, \label{eq60a}\\
			&&\dot{y}=1-by-x^{2},  \label{eq60b}\\
			&&\dot{z}=-x-cz. \label{eq60c}
			\end{eqnarray}
		\end{subequations}
Refs. \cite{chen2008dynamics}, \cite{chen2008nonlinear}, \cite{gao2009chaos}, \cite{jun2001study1}, \cite{jun2001study2} said that this financial dynamic model is formed of four sub-blocks: production, money, stock and labor force, and the dynamics  can be expressed as three first order differential equations. Three state variables of the system denote interest rate $x$, the investment demand $y$, and the price index $z$. To mention the constants, $a\geq0$ is the saving amount, $b\geq0$ is the cost per investment, and $c\geq0$ is the elasticity of demand of commercial markets. Two factors cause the major changes in the interest rate $x$: one of them is contradictions from the investment market, which is the surplus between investment and savings, and the other one is structural adjustment from the prices. This is expressed in \eqref{eq60a}. The rate of change of  $y$ is  related with the cost of investment and the interest rate as given in \eqref{eq60b}.  Change in $z$ is affected by inflation rates, therefore, at the same time,  it can be expressed  by the nominal interest rate and real interest rate, which is formulated in \eqref{eq60c} \cite{jun2001study1}.

The model we will be interested in is based on two existing models. When we focus on Ref. \cite{gao2009chaos}, we see that besides exhibiting the chaotic behaviour of the model \eqref{eq60} for some ranges of the parameters, by calculating the Lyapunov exponents; they also consider the case where there is a time delay feedback in the investment demand;
\begin{subequations}\label{main1}
	\begin{eqnarray}
&&\dot{x}=z(t)+[y(t)-a]x(t), \\
&&\dot{y}=1-by(t)-x^{2}+K[y(t)-y(t-\tau)], \\
&&\dot{z}=-x(t)-cz(t).
	\end{eqnarray}
\end{subequations}
Here $\tau \geq 0$ is the time delay and $K$ stands for the strength of the feedback.  Depending on the parameters, the system may have one or three equilibrium points. They perform the stability analysis of the system at one of the equilibrium points and occurrence of a Hopf bifurcation  in which the variable $y$ experiences periodic behavior is exhibited. If the constants $a,b,c$ of the system satisfy certain conditions, the authors show that the system is stable for $\tau\in [0,\tau_0)$, where $\tau_0$ is a critical value of the time delay, and the system undergoes a Hopf bifurcation when $\tau=\tau_0$.

The other model that we build our main problem upon is the hyperchaotic system of Ref.  \cite{yu2012dynamic}, which they formulate as
\begin{subequations}\label{main2}
\begin{eqnarray}
&&\dot{x}=z+(y-a)x+u, \\
&&\dot{y}=1-by-x^{2}, \\
&&\dot{z}=-x-cz,      \\
&&\dot{u}=-dxy-ku.
\end{eqnarray}
\end{subequations}
Basically, the authors of \cite{yu2012dynamic} state that the factors related to interest rate are relevant  not only to investment demand and price index but also to  the average profit margin: average profit margin and interest rate are proportional. By adding average profit margin as a new state variable $u$ to the system \eqref{eq60}, they obtain the system \eqref{main2}. This newly constructed system has an interesting property: While the system \eqref{eq60} has one positive  Lyapunov exponent for some range of the parameters, a sign for intrinsic chaotic behaviour, \eqref{main2} is shown in \cite{yu2012dynamic} to possess two positive Lyapunov  exponents for some region of the parameter space, which is mentioned in literature as a signal to hyperchaotic behaviour.

Motivated by the two works above, we consider the following system,
\begin{subequations}\label{main}
\begin{eqnarray}
			&&\dot{x}=z(t)+[y(t)-a]x(t)+u(t),      \\
			&&\dot{y}=1-by(t)-x^{2}(t)+K[y(t)-y(t-\tau)], \\
			&&\dot{z}=-x(t)-cz(t),   \\
			&&\dot{u}=-dx(t)y(t)-ku(t),
\end{eqnarray}
\end{subequations}
which is a combination of \eqref{main1} and \eqref{main2}, taking into account a time-delayed feedback in the investment demand variable $y$
and the effect of average profit margin simultaneously in \eqref{eq60}.
Here $K$ is the feedback strength, $K,\tau \geq 0$, and also $a,b,c,d,k$ are the nonnegative parameters of the system \eqref{main}.

To our knowledge, the existing literature does not consider the system \eqref{main} and in our analysis, we would like to answer the following questions:

\begin{itemize}
\item How is the qualitative behaviour of the system \eqref{main} around its fixed points?
\item When we follow the route in  \cite{gao2009chaos} and do the stability analysis of \eqref{main}, does the system undergo a Hopf bifurcation?
\item Can we analytically determine the critical value of $\tau_0$ that gives the bifurcation?
\item How is the effect of addition of the delay term on the stability of the system \eqref{main2}?
\end{itemize}

Therefore, the main purpose of this work is to search the dynamics of the financial model \eqref{main} by taking into account the effects of the additional variable and delay-feedback terms in \eqref{eq60}. After performing  stability analysis of the  constructed finance system, we theoretically demonstrate that the system undergoes a  Hopf bifurcation and this phenomenon is supported by  numerical simulations.

\noindent
The article  is organized as follows. In the following section we present a literature survey. Section 3 contains the main work, presenting  the stability analysis and the investigation of a Hopf bifurcation for  the constructed finance system at the fixed points. Bifurcation analyses are demonstrated by numerical simulations in Section 4. We conclude by some final remarks and future discussions.

\section{Related literature}\label{literaturereview}

The motivation of  this work is based on the system \eqref{eq60}, and, as we explained above, our model is a combination of \eqref{main1} and \eqref{main2}. In addition to the references mentioned before, in this subsection we will present a brief literature survey on these type of systems.

Refs. \cite {jun2001study1},  \cite{jun2001study2}  and \cite{junhai2008hopf}  consider the topological structure, Hopf bifurcation,  chaotic behaviour with different parameter combinations and the results of any change of the parameters on economy in  the  system \eqref{eq60}. Another study \cite{chen2008nonlinear} tackles with  \eqref{eq60} in view of fractional nonlinear models and its aim is to consider the chaotic behavior in fractional financial systems. Also,  Ref. \cite{zhao2011synchronization} considers synchronization strategies of a three-dimensional chaotic finance system.

By doing the shift $y\rightarrow y-\frac{1}{b}$ in the equation \eqref{eq60}, and adding a delay term to the first equation of the system,
\begin{subequations}
\begin{eqnarray}
&&\dot{x}=(\frac{1}{b}-a)x+z+xy+k(x(t-\tau)-x(t)), \\
&&\dot{y}=-by-x^{2}, \\
&&\dot{z}=-x-cz
\end{eqnarray}
\end{subequations}
is obtained, which is analysed in  \cite{kai2016chaotic}. Another version of this system appears in  \cite{yang2016bifurcation} as
\begin{subequations}
	\begin{eqnarray}
&&\dot{x}=-a(x+y)+K(x(t)-x(t-\tau)), \\
&&\dot{y}=-y-axz,  \\
&&\dot{z}=b+axy.
\end{eqnarray}
\end{subequations}

In Refs. \cite{ding2012hopf} and \cite{gao2009chaos}, time delay is added to the second equation of  \eqref{eq60}, the system becoming
\begin{subequations}
	\begin{eqnarray}
&&\dot{x}=z(t)+[y(t)-a]x(t), \\
&&\dot{y}=1-by(t)-x^{2}+K[y(t)-y(t-\tau)], \\
&&\dot{z}=-x(t)-cz(t),
	\end{eqnarray}
\end{subequations}
 to investigate the influence of the  time delay on investment demand $y$. Chen's system \cite{chen1999yet} is expressed as the following
\begin{subequations}
	\begin{eqnarray}
	&&\dot{x}=a(y-x), \\
	&&\dot{y}=(c-a)x-xz+cy, \\
	&&\dot{z}=xy-bz
	\end{eqnarray}
\end{subequations}
and by adding to the second equation of Chen system a time-delay term,   \cite{song2004bifurcation} obtains
\begin{subequations}
	\begin{eqnarray}
	&&\dot{x}=a(y-x), \\
	&&\dot{y}=(c-a)x-xz+cy+K(y(t)-y(t-\tau)),  \\
	&&\dot{z}=xy-bz
	\end{eqnarray}
\end{subequations}
and they study  the effect of the delayed feedback on Chen's system   and the existence of a Hopf bifurcation.

Another delayed financial  model is handled as follows in \cite{zhang2019hopf},
\begin{subequations}
	\begin{eqnarray}
	&&\dot{x}=(y-a)x+z(t-\tau),  \\
&&\dot{y}=1-by-x^{2},   \\
	&&\dot{z}=-x-cz
	\end{eqnarray}
	\end{subequations}
where $\tau$ represents price change delay.

In the studies \cite{chen2008dynamics} and  \cite{son2011delayed}, the authors construct the delayed financial system as follows
\begin{subequations}
	\begin{eqnarray}
	&&\dot{x}=z+(y-a)x+k_{1}\{x-x(t-\tau_{1})\},  \\
	&&\dot{y}=1-by-x^{2}+k_{2}\{y-y(t-\tau_{2})\},  \\
	&&\dot{z}=-x-cz+k_{3}\{z-z(t-\tau_{3})\}
	\end{eqnarray}
\end{subequations}
where $\tau_{1}, \tau_{2}$, and $\tau_{3}$ are time delays and $k_{1},k_{2}$, and $k_{3}$ demonstrate the strengths of the feedbacks. The aim is to investigate the  effect of delayed feedbacks on the financial system with time delay terms on the  interest rate, the investment demand and the price index of the financial system.

Another system  is constructed   in \cite{wu2010hopf} by adding a fourth variable $\omega$ to the system  Qi proposed in \cite{qi2005analysis}
\begin{subequations}
	\begin{eqnarray}
	&&\dot{x}=a(y-x)+eyz-k\omega \\
	&&\dot{y}=cx-dy-xz           \\
	&&\dot{z}=xy-bz               \\
	&&\dot{\omega}=rx+fyz
	\end{eqnarray}
\end{subequations}
 and the new system has chaotic or hyperchaotic behavior for wide frequency bandwith for suitable parameters.

Based on this literature, for the delayed model \eqref{main} we perform a stability analysis and exhibit  the existence of a Hopf bifurcation.
We provide several conditions for the stability of the system and analytically determine the critical delay value at which the system undergoes the bifurcation. Following this, by numerical simulations we show that the algebraic conditions on the coefficients of the system for the bifurcation to occur are indeed possible for a set of parameters.

\section{Stability of equilibrium points and Hopf bifurcation}
This section will include the stability analysis of our model \eqref{main}.
By solving the equations
\begin{subequations}\label{eq25}
	\begin{eqnarray}
 z^*+[y^*-a]x^*+u^*&=&0,    \\
 1-by^*-x^{*2}+K[y^*(t)-y^*(t-\tau)]&=&0,    \\
 -x^*-cz^*&=&0,   \\
 -dx^*y^*-ku^*&=&0,
\end{eqnarray}
\end{subequations}
and taking into account that $y^*=y^*(t)=y^*(t-\tau)$, we find the equilibrium points of the system.
The results are given in the following Lemmas.
\begin{lemma}
 In the situation
\begin{equation}
\frac{kb+abck+cd-ck}{c(d-k)} \leq0 ,
\end{equation}
the  system has a unique equilibrium,
\begin{equation}
 P_{0}\left(0,\frac{1}{b},0,0\right).
\end{equation}
\end{lemma}

\begin{lemma}
If the parameters of the system satisfy
\begin{equation}
\frac{kb+abck+cd-ck}{c(d-k)}>0 ,
\end{equation}
the  system has three equilibrium points;
\begin{align}
P_{0}\left(0,\frac{1}{b},0,0\right)
\end{align}
and
\begin{equation}
P_{1,2}\Bigg(\mp \sqrt{\frac{kb(1+ac)}{c(d-k)}+1}, \frac{(1+ca)k}{c(k-d)},
\pm  \frac{1}{c} \sqrt{\frac{kb(1+ac)}{c(d-k)}+1}, \mp \frac{d(1+ca)}{c(d-k)}\sqrt{\frac{kb(1+ac)}{c(d-k)}+1}\;\Bigg).
\end{equation}
\end{lemma}
\noindent
By the change of variables
\begin{equation}
X=x, \quad Y= y-\frac{1}{b}, \quad Z=z, \quad U=u
\end{equation}
the  equilibrium point $P_0$ is shifted to
\begin{equation}
P_{0}\left(0,0,0,0\right).
\end{equation}
\noindent
After this transformation, the new system can be arranged as;
\begin{subequations}\label{sysX}
\begin{eqnarray}
\dot{X}&=&Z(t)+[Y(t)+\frac{1}{b}-a]X(t)+U(t),      \\
\dot{Y}&=&-bY(t)-X^{2}(t)+K[Y(t)-Y(t-\tau)],  \\
\dot{Z}&=&-X(t)-cZ(t)   \\
\dot{U}&=&-dX(t)[Y(t)+\frac{1}{b}]-kU(t).
\end{eqnarray}
\end{subequations}
\subsection{Stability analysis and Hopf bifurcation for $P_{0}$}

We work with the system  \eqref{sysX}, for which the equilibrium point $P_0$ is the origin.
Remember that the characteristic equation for a delayed system is
\begin{equation}\label{jacob}
|J_{0}+e^{-\lambda\tau} J_{\tau} - \lambda I|=0.
\end{equation}
 At the equilibrium point $P_0(0,0,0,0)$, we find $J_0$ and $J_\tau$ for Eq. \eqref{sysX} to be
\begin{equation}
J_{0,0} =\left( \begin{array}{cccc}
\frac{1}{b}-a & 0 & 1 &  1 \\
0 & -b+K & 0 & 0 \\
-1 & 0 & -c & 0\\
-\frac{d}{b} & 0 & 0 & -k \\
\end{array} \right)
\end{equation}
and
\begin{equation}
J_{\tau,0} =
\left( \begin{array}{cccc}
0 & 0 & 0 &  0 \\
0 & -K & 0 & 0 \\
0 & 0 & 0 & 0  \\
0 & 0 & 0 & 0  \\
\end{array} \right).
\end{equation}
\noindent Solving \eqref{jacob}, the characteristic equation is obtained, which appears to be a fourth degree exponential polynomial equation as follows
\begin{equation}\label{char}
\Big[\lambda-(-b+K-Ke^{-\lambda\tau})\Big](\lambda^{3}+p_{1}\lambda^{2}+p_{2}\lambda+p_{3})=0
\end{equation}
where
\begin{subequations}\label{RH}
\begin{eqnarray}
\ p_{1}&=&k+a+c-\frac{1}{b}, \\
\ p_{2}&=&ck+ak+ac-\frac{k+c-d}{b}+1, \\
\ p_{3}&=&(1+ac-\frac{c}{b})k+\frac{cd}{b}.
\end{eqnarray}
\end{subequations}
\begin{remark}
Let us note that the characteristic equation of the system \eqref{main1} appearing in \cite{gao2009chaos} is
\begin{equation}\label{chargao}
\Big[\lambda-(-b+K-Ke^{-\lambda\tau})\Big](\lambda^{2}+P\lambda+Q)=0.
\end{equation}
We see that \eqref{char} and \eqref{chargao} share the same transcendental part, whereas in our case, the algebraic part is of third degree, and more complicated than that of  \cite{gao2009chaos}
\end{remark}

In the following result we are considering the algebraic part of the characteristic equation \eqref{char}.
\begin{lemma}
According to the Routh-Hurwitz criterion, when the conditions
\begin{equation}\label{HRW}
p_1>0, \quad p_3 > 0, \quad p_1 p_2 > p_3
\end{equation}
hold, the three roots of the characteristic equation \eqref{char} originating from the algebraic term in the second parenthesis, have negative real parts; i.e., the roots are on the left half plane, for all $\tau\geq 0$.
\end{lemma}
\noindent
\begin{remark}
When $\tau=0$,  the  transcendental part of \eqref{char} reduces to
\begin{equation}\label{eq28}
\lambda+b=0.
\end{equation}
$b$ represents the cost per investment.  When $b>0$, the root of  the equation \eqref{eq28} is negative. Therefore, in case $\tau=0$, when the conditions of Lemma 3  are met,  the system  is stable at the equilibrium $P_{0}$.
\end{remark}
See that Lemma 3 guarantees that when the conditions in \eqref{HRW} are satisfied, for all values of $\tau \geq 0$, the eigenvalues originating from the algebraic part of \eqref{char} have negative real parts, therefore stability is at our disposal for this part of the characteristic equation. Then we just need to analyze the root of the transcendental part of \eqref{char},  which is as follows
\begin{equation}\label{eq29}
\lambda+b-K+Ke^{-\lambda\tau}=0.
\end{equation}
\noindent
The machinery below in this subsection shares the same calculations and lines with the computations in \cite{gao2009chaos}. We repeat these calculations here for completeness. Assume that the root of the transcendental equation is of the form
\begin{equation}
\lambda(\tau)=\alpha(\tau)+\beta(\tau)i.
\end{equation}
Suppose that for some critical value $\tau=\tau^{*}$, we have  $\lambda(\tau^*)=0$ and  $\beta(\tau^*)\neq0$. Then the system
\eqref{sysX} undergoes a Hopf bifurcation at the origin provided the transversality condition is satisfied. It is clear that, if $\lambda=i\omega$ ($\omega>0$) is a root of \eqref{eq29},  it must satisfy
\begin{equation}
\ i\omega+b-K+K(\cos\omega\tau-i\sin\omega\tau)=0.
\end{equation}
\noindent
If the imaginary and real parts are separated, we obtain;
\begin{subequations}\label{eq30}
	\begin{eqnarray}
\omega & = & K\sin\omega\tau, \\
\ K-b & = & K \cos\omega\tau.
\end{eqnarray}
\end{subequations}
Eliminating the trigonometric terms, one obtains
\begin{equation}
\omega^{2}=2Kb-b^{2}
\end{equation}
and hence
\begin{equation}\label{eq31}
\omega_{+}=\sqrt{2Kb-b^{2}}.
\end{equation}
It is clear that if $K> b/2 $,\; $\omega>0$ is determined  uniquely, and if $K\leq\ b/2 $, \; there is no such $\omega$.
Using \eqref{eq30} we get
\begin{equation}\label{tauj}
\tau_{j}=\frac{1}{\omega_{+}}\arccos\frac{K-b}{K}+\frac{2j\pi}{\omega_{+}}, \quad j=0,1,2,\dots
\end{equation}
\noindent
We have proven that $\tau=\tau_j$, $\lambda(\tau_j)=i\omega_+$ is a pure imaginary root of the transcendental equation \eqref{eq29}.

In order to show the existence of the Hopf bifurcation, we need to check the transversality condition. This is done in the following Lemma, which we take from Ref. \cite{gao2009chaos}.
\begin{lemma}
$\lambda(\tau_j)$ satisfies the transversality condition; that is,   $ \displaystyle \frac{d\mathrm{Re}(\lambda(\tau_{j}))}{d\tau}>0$ \,for  $j=0,1,2,\ldots$.
\end{lemma}
\textit{Proof}: Considering $\lambda=\lambda(\tau)$ in \eqref{eq29},
\begin{equation}\label{eq33}
\lambda+b-K+Ke^{-\lambda\tau}=0,
\end{equation}
and evaluating the derivative of the equation with respect to $\tau$, one obtains
\begin{equation}
\frac{d\mathrm{Re}(\lambda)}{d\tau} \Big\vert_{\tau=\tau_{j}}  =\frac{\omega_{+}^{2}}{(\cos \omega_{+} \tau_{j}-K \tau_{j})^{2}+(\sin \omega_{+} \tau_{j})^{2}}>0;
\end{equation}
which proves the Lemma. For details of the calculation, see Ref. \cite{gao2009chaos}.
$ \qed $

\noindent
To analyze the roots of the exponential polynomial equation \eqref{eq29}, the result of the following Lemma which is proved by Ruan \& Wei \cite{ruan2003zeros} is needed.
\begin{lemma}\label{ruan}
Consider the exponential polynomial equation
\begin{eqnarray}
    P(\lambda,  e^{-\lambda \tau_{1}},\ldots, e^{-\lambda \tau_{m}})
    &=&\lambda^{n}+p_{1}^{(0)} \lambda^{n-1}+\ldots+p_{n-1}^{(0)} \lambda+p_{n}^{(0)}
    +[p_{1}^{(1)} \lambda^{n-1}+\ldots +p_{n-1}^{(1)} \lambda +p_{n}^{(1)}]e^{-\lambda \tau_{1}}   \nonumber \\
    &+&\ldots +[p_{1}^{(m)} \lambda^{n-1}+\ldots+p_{n-1}^{(m)} \lambda +p_{n}^{(m)}]e^{-\lambda \tau_{m}}=0,
\end{eqnarray}
where $\tau_{i} \geq 0$ $(i=1,2,...,m)$ and $p_{j}^{(i)}$ $(i=0,1,2,...,m; j=1,2,...,n)$ are constants. As $(\tau_{1},\tau_{2},...,\tau_{m})$ vary, the sum of the order of the zeros of $P(\lambda, e^{-\lambda \tau_{1}},...,e^{-\lambda \tau_{m}})$ on the open right half-plane can change only if a zero appears on or crosses the imaginary axis.
\end{lemma}
\begin{remark}\label{rmk2}
Let us summarize the steps we have gone through so far and comment on how to arrive at Theorem 1.
\end{remark}
\begin{itemize}
    \item When $\tau=0$: If condition \eqref{HRW} of Lemma 3 hold and if $b>0$, then all the eigenvalues of the linearization of the system \eqref{sysX} has negative real parts, hence \eqref{sysX} is stable at $P_0$ when $\tau = 0$. When $\tau = 0$, \eqref{char} is an algebraic equation of degree four. Since all the eigenvalues have negative real parts,  the sum of the multiplicities of zeros of \eqref{char} (let us call this $LS$) on the left half plane is equal to $4$, write $LS=4$.  Since there is no eigenvalue with a positive real part, the sum of multiplicities of the eigenvalues on the right half plane (let us call this $RS$) is $0$, write $RS=0$.
    \item Lemma 5 is an only if statement, and it says, if $RS$ changes, then there arises a pure imaginary root of \eqref{char}. This is equivalent to the following statement: If \eqref{char} does not have any imaginary root, then $RS$ does not change.
    \item Suppose $K>b/2$. Then, we have a pure imaginary root of the characteristic equation.
    \begin{itemize}
    \item The smallest value of $\tau$ such that \eqref{char} has a pure imaginary root is $\tau=\tau_0$, where $\lambda(\tau_0)=i \omega_+$. Therefore, when $\tau\in[0,\tau_0)$,  $RS$ does not change and remains the same as $RS=0$. Hence, when  $\tau\in[0,\tau_0)$, there is no eigenvalue with positive real part. $\lambda = 0$ is already not an eigenvalue  since we put the conditions $b>0$ and $p_3>0$. Hence when $\tau\in[0,\tau_0)$, all the eigenvalues have negative real parts, and the system is stable at the equilibrium point $P_0$.
    \item We have shown that $\mathrm{Re}\big(\lambda(\tau_0)\big)=\mathrm{Re}(i\omega_+)=0$ \quad and \quad $\displaystyle \frac{d\mathrm{Re}(\lambda(\tau))}{d\tau}|_{\tau=\tau_0}>0$. This means, at $\tau=\tau_0$,  $\mathrm{Re}(\lambda(\tau))$ is an increasing function of  $\tau$. $\mathrm{Re\big(\lambda(\tau)\big)}$ must pass from negative to positive values at $\tau=\tau_0$. Hence, on the interval $(\tau_0,\tau_1)$, $\mathrm{Re(\lambda)}$ takes on positive values. For some value of $\tau$  in $(\tau_0,\tau_1)$, we have $RS\geq 1$. Since on the interval $(\tau_0,\tau_1)$  there does not occur a pure imaginary root, $RS$ does not change, hence $RS\geq 1$  for every $\tau \in (\tau_0,\tau_1)$.
    Therefore, for $\tau \in (\tau_0,\tau_1)$ the system has at least one eigenvalue with positive real part and it is unstable at the equilibrium point $P_0$.
    \end{itemize}
    \item Suppose $K\leq b/2$. Then, we do not have a pure imaginary root of the characteristic equation.
    \begin{itemize}
        \item If condition \eqref{HRW} of Lemma 3 hold and if $b>0$, since no imaginary root occurs for any value of $\tau\geq 0$, there does not appear any eigenvalue on the right half plane and the system remains stable for all values of $\tau\geq 0$.
    \end{itemize}
\end{itemize}

Based on the arguments that we tried to explain above, considering the Lemmas 1-5, the following Theorem can be obtained.
\begin{theorem} We assume that the conditions of \eqref{HRW}  hold and $b>0$.

\vspace{0.3cm}
\noindent If $K>b/2$,
\begin{itemize}
\item[(i)] The equilibrium point $P_{0}$  of the system \eqref{sysX} is stable for $\tau\in [0,\tau_0)$ while it behaves unstable for $\tau \in (\tau_0,\tau_1)$.
\item[(ii)] The system \eqref{sysX} undergoes a Hopf bifurcation at the equilibrium point $P_{0}$ when $\tau=\tau_0$.
\end{itemize}
\vspace{0.3cm}
\noindent If $K\leq b/2$,
\begin{itemize}
\item[(iii)] The equilibrium point $P_{0}$  of the system \eqref{sysX} is stable for  $\tau \geq 0$.
\end{itemize}
\end{theorem}
\subsection{Stability analysis and Hopf bifurcation for $P_{1}$}
In the previous subsection, the stability condition of the system \eqref{main} is considered at the equilibrium point $P_{0}={(0,\frac{1}{b},0,0)}$, or equivalently, of the system \eqref{sysX} at the point $P_0(0,0,0,0)$. There remains to consider the stability of the other two equilibrium points. We shall deal with $P_1$ only, the analysis follows similar lines for the equilibrium point $P_2$.

Remember that the equilibrium point $P_1$ was
\begin{equation}
P_{1}\Bigg(\theta, \; \frac{k(1+ac)}{c(k-d)}, -\frac{1}{c} \theta, \,  \frac{d(1+ac)}{c(d-k)}\theta\;\Bigg),
\end{equation}
where
\begin{equation}
\theta=\sqrt{\frac{kb(1+ac)}{c(d-k)}+1}.
\end{equation}
The related Jacobians $J_0$ and $J_\tau$ at the  equilibrium $P_{1}$  are found to be
 \begin{equation}
 J_{0,1} =
 \left( \begin{array}{cccc}
 \frac{k+acd}{c(k-d)} & \theta & 1 &  1 \\
 -2\theta &  -b+K & 0 & 0 \\
 -1 & 0 & -c & 0\\
 -\frac{dk(1+ac)}{c(k-d)} & -d\theta & 0 & -k \\
 \end{array} \right)
 \end{equation}
 and
 \begin{equation}
 J_{\tau,1} =
 \left( \begin{array}{cccc}
0 & 0 & 0 &  0 \\
0 & -K & 0 & 0 \\
0 & 0 & 0 & 0  \\
0 & 0 & 0 & 0  \\
\end{array} \right).
 \end{equation}
The characteristic equation
\begin{align}
|J_{0}+e^{-\lambda\tau} J_{\tau} - \lambda I|=0
\end{align}
yields
\begin{equation}\label{sys3}
\lambda^{4}+a_{1}\lambda^{3}+b_{1}\lambda^{2}+c_{1}\lambda+d_{1}+(a_{2}\lambda^{3}+b_{2}\lambda^{2}+c_{2}\lambda)e^{-\lambda\tau}=0
\end{equation}
where
\begin{align}
a_{1}&=b+c+k-K+\frac{acd+k}{c(d-k)}, \nonumber\\
b_{1}&=1+ck+2\theta^2+(c+k)(b-K)+\frac{1}{c (d-k)}\Big[acd(b+c-K)+k(b+c-d+k-K)\Big], \nonumber\\
c_{1}&=\frac{b-K}{c(d-k)}\Big[cd+k(k-d)+c^2\big(ad+k(d-k)\big)\Big]+2(c-d+k)\theta^{2},  \nonumber\\
d_{1}&=2c(k-d)\theta^{2},	\\
a_{2}&=K,  \nonumber\\
b_{2}&=\Big[c+k+\frac{acd+k}{c(d-k)}\Big]K, \nonumber\\
c_{2}&=\Big[1+ck+\frac{adc^2+k(c-d+k)}{c(d-k)}\Big]K. \nonumber
\end{align}

\noindent
We see that the characteristic equation is different than that of $P_0$, and we will make use of works of  Ruan and Wei \cite{ruan2001zeros},  \cite{ruan2003zeros} and Li and Wei \cite{li2005zeros} in order to discover the distribution of zeros of Eq. \eqref{sys3}, which is  a fourth degree transcendental polynomial equation.

\noindent
If $i\om$, $(\om>0)$ is a root of   Eq. \eqref{sys3},  it must satisfy
\begin{equation}\label{eq80}
\om^{4}-a_{1}\om^{3}i-b_{1}\om^{2}+c_{1}\om i +d_{1} +(-a_{2}\om^{3}i-b_{2}\om^{2}+c_{2}\om i)e^{-i\om \tau}=0.
\end{equation}
After separating the real and imaginary parts of the equation, we get
\begin{align}
\om^{4}-b_{1}\om^{2}+d_{1}&=(a_{2}\om^{3}-c_2 \om)\sin(\om \tau)+b_{2}\om^{2}\cos(\om \tau),   \label{sincos1}\\
a_{1} \om^{3}-c_{1}\om&=(c_2\om-a_{2}\om^{3})\cos(\om \tau)+b_{2}\om^{2}\sin(\om \tau).    \label{sincos2}
\end{align}
Taking the squares of both equations and adding up obtain
\begin{equation}\label{eq37}
\om^{8}+(a_{1}^{2}-2b_{1}-a_{2}^{2})\om^{6}+(b_{1}^{2}+2d_{1}-2a_{1}c_{1}-b_{2}^{2}+2a_{2}c_{2})
\om^{4}+(c_{1}^{2}-2b_{1}d_{1}-c_{2}^{2})\om^{2}+d_{1}^{2}=0.
\end{equation}
Let $z=\om^{2}$ and denote $p=a_{1}^{2}-2b_{1}-a_{2}^{2}$ , $q=b_{1}^{2}+2d_{1}-2a_{1}c_{1}-b_{2}^{2}+2a_{2}c_{2}$,
$u=c_{1}^{2}-2b_{1}d_{1}-c_{2}^{2}$ and $v=d_{1}^{2}$. Then Eq. \eqref{eq37} becomes
\begin{equation}\label{eq40}
z^{4}+pz^{3}+qz^{2}+uz+v=0.
\end{equation}
\begin{remark}
The transcendental-polynomial characteristic equation \eqref{sys3}  is different and more complicated than the
characteristic equation considered in \cite{li2005zeros} .  Their characteristic equation is in the form
\begin{equation}\label{weichar}
\lambda^4+a\lambda^3+b\lambda^2+c\lambda+d+re^{-\lambda \tau}=0,
\end{equation}
which appears as Eq. (2.1) in \cite{li2005zeros}. Although our characteristic equation \eqref{sys3} was different than theirs, we followed their lines
in the search of a root $\lambda = i \om$, $\om > 0 $ and obtained exactly the same fourth-degree equation \eqref{eq40} up to a difference in all of the constants, of course. Therefore, apart from this line, we will adapt the development starting by Eq. (2.5) of   \cite{li2005zeros}.
\end{remark}
\noindent
Let us call
\begin{equation}\label{eq48}
h(z)=z^{4}+pz^{3}+qz^{2}+uz+v.
\end{equation}
\noindent
After differentiating $h(z)$, we have $h'(z)=4z^{3}+3pz^{2}+2qz+u$. Set
\begin{equation}\label{eq39}
4z^{3}+3pz^{2}+2qz+u=0.
\end{equation}
Let $y=z+\frac{p}{4}$. Then equation \eqref{eq39} becomes
\begin{equation}\label{eq41}
y^{3}+p_{1}y+q_{1}=0,
\end{equation}
where
$\displaystyle p_{1}=\frac{q}{2}-\frac{3}{16}p^{2}$, \,  $\displaystyle  q_{1}=\frac{p^{3}}{32}-\frac{pq}{8}+\frac{u}{4}$.

\noindent Let us define
\begin{align}
&D=(\frac{q_{1}}{2})^{2}+ (\frac{p_{1}}{3})^{3}, \quad  \sigma=\frac{-1+\sqrt{3}i}{2},\\
&y_{1}=\sqrt[3]{-\frac{q_{1}}{2}+\sqrt{D}}+\sqrt[3]{-\frac{q_{1}}{2}-\sqrt{D}}, \\
&y_{2}=\sqrt[3]{-\frac{q_{1}}{2}+\sqrt{D}\sigma}+
\sqrt[3]{-\frac{q_{1}}{2}-\sqrt{D}\sigma^{2}}, \\
&y_{3}=\sqrt[3]{-\frac{q_{1}}{2}+\sqrt{D}\sigma^{2}}+
\sqrt[3]{-\frac{q_{1}}{2}-\sqrt{D}\sigma}, \\
&z_{i}=y_{i}-\frac{p}{4}, \, i=1,2,3.
\end{align}

\begin{lemma}
Observe that $v=d_1^2\geq0$.
\begin{itemize}
\item[(i)]  When $D\geq0$, the equation \eqref{eq40} has positive roots  iff $z_{1}>0$ and $h(z_{1})<0$;
\item[(ii)] When $D<0$, the equation \eqref{eq40} has positive root iff there exists at least one $z^{*}\in{z_{1},z_{2},z_{3}}$, such that $z^{*}>0$ and $h(z^{*})\leq0$. (Adapted from \cite{li2005zeros})
\end{itemize}
\end{lemma}

\noindent
\textit{Proof}: Available in Ref.  \cite{li2005zeros}. The statement of the Lemma 6 for our case is slightly different than theirs. $ \qed $

\noindent
Assume that Eq. \eqref{eq40} has positive roots. We can suppose that it has four positive roots which we will note as $z_{i}^{*}$, $i=1,2,3,4$. It follows that the equation \eqref{eq37} has four positive roots, denoted by  $\om_{i}=\sqrt{z_{i}^{*}}$, $i=1,2,3,4$.
From \eqref{sincos1} and \eqref{sincos2} we solve
\begin{align}\label{taucr2}
\tau_{k}^{(j)}=\frac{1}{\om_{k}}\Big\{\arccos \Big[& \big((b_2-a_1a_2)w_k^4 + (a_1c_2+a_2c_1-b_1b_2)w_k^2+b_2d_1 \\
                                &-c_1c_2  \big)  / \big( a_2^2w_k^4+(b_2^2-2a_2c_2)w_k^2+c_2^2  \big)\Big] +2(j-1)\pi\Big\} \nonumber
\end{align}
for $k=1,2,3,4$; \; $j=1,2,\ldots $. Then $\mp i\om_{k}$ is a pair of purely imaginary roots of equation \eqref{sys3} when $\tau=\tau_{k}^{(j)}$, $k=1,2,3,4; \, j=1,2,...$. Obviously,
\begin{equation}
\lim_{j\to\infty}\tau_{k}^{(j)}=\infty, \qquad  k=1,2,3,4.
\end{equation}
Then, it can be defined that
\begin{equation}\label{eq62}
\tau_{0}=\tau_{k_{0}}^{(j_{0})}=\mathrm{min}_{_{1\leq k \leq 4},_{1 \leq j}}{\{\tau_{k}^{(j)}}\},\quad \om_{0}=\om_{k_{0}}, \quad z_{0}=z_{k_{0}}^{*}.
\end{equation}
\begin{lemma}\label{eq75}
 Assume that $a_{1}+a_{2}>0$,\, $(a_{1}+a_{2})(b_{1}+b_{2})-(c_{1}+c_{2})>0$,\, $d_{1}>0$ \, and $(c_{1}+c_{2})[(a_{1}+a_{2})(b_{1}+b_{2})-(c_{1}+c_{2})]-(a_{1}+a_{2})^{2}d_{1}>0$.
\end{lemma}
\begin{itemize}
	\item [(i)] For $\tau\in[0,\tau_{0})$, if one of following conditions holds: \emph{(a)} $D\geq0$, $z_{1}>0$ and $h(z_{1})\leq0$; \emph{(b)} $D<0$ and the roots of the equation $z^{*}\in\{ z_{1},z_{2},z_{3}\}$ are such that $\exists\, z^{*}>0$ and $h(z^{*})\leq0$, then all roots of equation \eqref{sys3} have negative real parts.
	\item [(ii)]  For all $\tau\geq0$, if the conditions $(a) \, and \,(b)$ of $(i)$ are not satisfied, then all roots of equation \eqref{sys3} have negative real parts.
	\end{itemize}

\textit{Proof}: When $\tau=0$, equation \eqref{sys3} becomes
\begin{equation}\label{eq45}
\lambda^{4}+(a_{1}+a_2)\lambda^{3}+(b_{1}+b_2)\lambda^{2}+(c_{1}+c_2)\lambda+d_{1}=0.
\end{equation}
According to the Routh-Hurwitz criterion, all roots of equation \eqref{eq45} have negative real  parts if and only if
\begin{align}
&a_{1}+a_{2}>0, \quad  d_{1}>0, \quad  (a_{1}+a_{2})(b_{1}+b_{2})-(c_{1}+c_{2})>0,   \nonumber \\
&(c_{1}+c_{2})[(a_{1}+a_{2})(b_{1}+b_{2})-(c_{1}+c_{2})]-(a_{1}+a_{2})^{2}d_{1}>0.
\end{align}
Now, under this condition, at $\tau = 0$, \eqref{sys3} has no eigenvalue with positive real part, $RS=0$.
Lemma 6 proposes if and only if conditions for \eqref{sys3} to have a  pure imaginary root. If $(a)-(b)$ are not satisfied, Lemma 6 says that \eqref{sys3} has no  pure imaginary root for all $\tau\geq0$. Therefore, according to Lemma 5, $RS$ does not change for \eqref{sys3} for all $\tau \geq 0$, which means it does not have any eigenvalue with positive real part for   $\tau \geq 0$. $\lambda = 0$ is not an eigenvalue since we assumed $d_1>0$. Therefore, all the eigenvalues are with negative real parts for $\tau \geq 0$. This proves \emph{(ii)}.

Now let us prove \emph{(i)}. Since one of \emph{(a)-(b)} holds, then by Lemma 6, there is a pure imaginary root of \eqref{sys3} at $\tau = \tau_0$. Since there is no pure imaginary root of \eqref{sys3} for all values of $\tau\in [0, \tau_0)$, by Lemma 5, $RS$ does not change on $0\leq \tau < \tau_0$, hence  all the eigenvalues of \eqref{sys3} are with negative real parts when $\tau\in [0, \tau_0)$. This proves \emph{(i)}. $ \qed $

Let
\begin{equation}\label{eq61}
\lambda(\tau)=\alpha(\tau)+i\om(\tau)
\end{equation}
be the root of equation \eqref{sys3} satisfying $\alpha(\tau_{0})=0$,\, $\om(\tau_{0})=\om_{0}$.
\begin{lemma}
Assume that $h'(z_{0})\neq0$. Then, $\mp i\omega_{0}$ is a simple (i.e., not multiple) pure imaginary root of the equation \eqref{sys3} when $\tau=\tau_{0}$.  Additionally, if the conditions of Lemma 7-(i) hold, the following transversality condition is satisfied:	
\begin{equation}
\frac{d(\mathrm{Re}\lambda(\tau))}{d\tau}\Big\vert_{\tau=\tau_{0}}\neq0
\end{equation}
and the sign of $d(\mathrm{Re}\lambda(\tau))/d\tau|_{\tau=\tau_{0}}$ is consistent with that of $h'(z_{0})$.
\end{lemma}
\noindent
\textit{Proof}: For the proof of this Lemma, we will use the results available in the proof of Lemma 7 of Ref. \cite{zhang2019hopf}.  Denote
\begin{align}
&R(\lambda)=\lambda^{4}+a_{1}\lambda^{3}+b_{1}\lambda^{2}+c_{1}\lambda+d_{1}, \\
&Q(\lambda)=a_{2}\lambda^{3}+b_{2}\lambda^{2}+c_{2}\lambda.
\end{align}
Then \eqref{sys3} can be represented as
\begin{equation}\label{eq50}
R(\lambda)+Q(\lambda)e^{-\lambda \tau}=0,
\end{equation}
and \eqref{eq37} can be written as the following:
\begin{equation}\label{eq51}
R(i\om)\bar{R}(i\om)-Q(i\om)\bar{Q}(i\om)=0.
\end{equation}
Then, together with \eqref{eq40} and \eqref{eq48}, we get
\begin{equation}\label{eq49}
h(\om^{2})=R(i\om)\bar{R}(i\om)-Q(i\om)\bar{Q}(i\om).
\end{equation}
Differentiating both sides of \eqref{eq49} with respect to $\om$, we obtain
\begin{equation}\label{eq52}
2\om h'(\om^{2})=-i\{R\bar{R'}-R'\bar{R}+	Q'\bar{Q}-Q\bar{Q'}\}.
\end{equation}

If $i\om_{0}$ is not a  simple root, i.e., if it is a multiple root,  then it must satisfy
\begin{equation}
\frac{d}{d\lambda}[R(\lambda)+Q(\lambda)e^{-\lambda \tau}]|_{\lambda=i\om_{0}}=0,
\end{equation}
that is,
\begin{equation}
R'(i\om_{0})+Q'(i\om_{0})e^{-i\om_{0}\tau_{0}}-\tau_{0}Q(i\om_{0})e^{-i\om_{0}\tau_{0}}=0.
\end{equation}
With \eqref{eq50}, we have
\begin{equation}
\tau_{0}=\frac{Q'(i\om_{0})}{Q(i\om_{0})}-\frac{R'(i\om_{0})}{R(i\om_{0})}.
\end{equation}
Thus, following the analysis in \cite{zhang2019hopf}, using \eqref{eq51} and \eqref{eq52}, one gets
\small
\begin{equation}\label{ref1}
\mathrm{Im}(\tau_{0})
=\frac{\om_{0}h'(\om_{0}^{2})}{|R(i\om_{0})|^{2}}.
\end{equation}
This yields $h'(\om_{0}^{2})=0$, since $\tau_{0}$ is real and $\mathrm{Im}(\tau_{0})=0$. We have a contradiction to the assumption $h'(\om_{0}^{2})\neq0$. This is the proof of  the first conclusion.

\noindent
Now, we need to prove that $\displaystyle\frac{d(\mathrm{Re}(\lambda\tau))}{d\tau}\Big\vert_{\tau=\tau_{0},\lambda=i\om_{0}}\neq0$.

\noindent
Differentiating both sides of  \eqref{eq50} with respect to $\tau$,
\begin{equation}
\frac{d}{d\tau}[R(\lambda)+Q(\lambda)e^{-\lambda\tau}]=0
\end{equation}
we obtain;
\begin{equation}
\frac{d\lambda}{d\tau}[R'(\lambda)+Q'(\lambda)e^{-\lambda\tau}-\tau Q(\lambda)e^{-\lambda\tau}]-\lambda Q(\lambda)e^{-\lambda\tau}=0
\end{equation}
which implies
\begin{equation}\label{ref2}
\frac{d\lambda}{d\tau}=\frac{\lambda[-R(\lambda)\bar{R'}(\lambda)e^{\lambda \tau}e^{\bar{\lambda}\tau}+Q(\lambda)\bar{Q'}(\lambda)-\tau |Q(\lambda)|^{2}]}{|R'(\lambda)e^{\lambda \tau}+ Q'(\lambda)-
\tau Q(\lambda)|^{2}}.
\end{equation}
It follows together with \eqref{eq52} that
\small
\begin{equation}\label{ref3}
\frac{d(\mathrm{Re}\lambda(\tau))}{d\tau}\Big\vert_{\tau=\tau_{0},\lambda=i\om_{0}}
=\frac{\om_{0}^{2}h'(\om_{0}^{2})}{|R'(i\omega_0)e^{i\omega_0 \tau_0}+ Q'(i\omega_0)-\tau_0 Q(i\omega_0)|^{2}}\neq0. \qquad \qed
\end{equation}
For the intermediate calculations leading to the results in equations \eqref{ref1}, \eqref{ref2}, \eqref{ref3} we refer the reader to
\cite{zhang2019hopf}. Now we can state the following Theorem.
\begin{theorem}
The values of  $\omega_{0}$, $\tau_{0}$, $z_{0}$ and $\lambda(\tau)$ are given in equations \eqref{taucr2}, \eqref{eq62} and \eqref{eq61}. Suppose that
\begin{itemize}
\item[$\bullet$] $a_{1}+a_{2}>0$, $(a_{1}+a_{2})(b_{1}+b_{2})-(c_{1}+c_{2})>0$,\quad $d_{1}>0$,
\item[$\bullet$] $(c_{1}+c_{2})[(a_{1}+a_{2})(b_{1}+b_{2})-(c_{1}+c_{2})]-(a_{1}+a_{2})^{2}d_{1}>0$.
\end{itemize}
\begin{itemize}
\item[(i)] When the assumptions (a)-(b) of Lemma \eqref{eq75} do not hold, the roots of equation \eqref{sys3} have negative real parts for all $\tau\geq0$. The system \eqref{main} is stable for all $\tau \geq 0$.
\item[(ii)] If   either (a) or (b) is satisfied, roots of  Eq. \eqref{sys3} have negative real parts for $\tau\in[0,\tau_{0})$. It is obtained  that one of  the roots of \eqref{sys3} is $\mp\omega_{0}i$ and other roots have negative real parts  when $\tau=\tau_{0}$ and $h'(z_0)\neq 0$. Additionally, $\displaystyle \frac{d\mathrm{Re}\lambda(\tau_{0})}{d\tau}>0$ and if  $\tau_{1}$ is taken the first value of $\tau>\tau_{0}$
 such that Eq. \eqref{sys3} has purely imaginary root, Eq. \eqref{sys3} has at least one root with  positive real part for $(\tau_0,\tau_1)$. Therefore, at the equilibrium point $P_1$, the system \eqref{main},
 \begin{itemize}
     \item is stable for $\tau \in [0,\tau_0)$,
     \item undergoes a Hopf bifurcation when $\tau=\tau_0$,
     \item is unstable for $\tau \in (\tau_0,\tau_1)$.
 \end{itemize}
\end{itemize}
\end{theorem}

Let us note that arriving at this Theorem 2 from Lemmas 5-8 is through the same lines of reasoning we tried to explain in Remark \ref{rmk2} for Thorem 1; we do not repeat them  here. In Lemma 8, it is shown that $\lambda = \mp i \omega_0$ is a simple pure imaginary root of the characteristic equation \eqref{sys3}. Being a simple root of the eigenvalue equation, not a multiple root, is a necessary condition for the existence of the derivative $\displaystyle \frac{d(\mathrm{Re}\lambda(\tau))}{d\tau}$, from which we check the transversality condition for a Hopf bifurcation to occur.
\section{Numerical simulations of the system}
Having obtained these theoretical results, we support our findings with plots of the time series of the dependent variables of the system for different values of the time-delay $\tau$. Our simulations are  in well accordance with the theoretical findings of the previous section. We give the dependent variable versus time plots for values of the time-delay $\tau$ for which
\begin{itemize}
    \item $\tau < \tau_0$, and the system is stable at the fixed points $P_0$, $P_1$,
    \item $\tau = \tau_0$, and the system undergoes a Hopf bifurcation exhibited by the periodic behaviour in the variable $y$ for  $P_0$ and in $x$ and $y$ for $P_1$,
    \item $\tau>\tau_0$, and the system becomes unstable, demonstrated by the unbounded development in the graph of the variable $y$ for  $P_0$ and in $x$ and $y$ for $P_1$.
\end{itemize}

We use the dde23 solver of MATLAB$^\circledR$ for the delay differential equation systems to solve the system \eqref{main} numerically.

In Figure 1, the  parameter values and initial conditions are chosen for testing the theoretical results which were obtained for the critical point $P_0(0,0,0,0)$ of the shifted system \eqref{sysX}. We consider the parameters $a=5$, $b=0.4$, $c=1.5$, $d=0.2$, $k=0.17$ and $K=1$ which satisfy the bifurcation criteria of Theorem 1. Also, the initial conditions are taken as $X(0)=1$, $Y(0)=2$, $Z(0)=U(0)=0.5$. The critical value of the delay is calculated from \eqref{tauj} as $\tau_0=1.15912$.  Numerical results given in (a), (b) and (c) are examples for  stable, Hopf bifurcation and unstable cases, respectively. It is also observed that the transition from stability to instability through the Hopf bifurcation is experienced only in the $y$ variable for the equilibrium point $P_0$.  These simulations are compatible with predictions of our theoretical results.

In the case of Figure 2, for the equilibrium point $P_1$ for the system \eqref{main}, we choose the parameters $a=0.2$, $b=0.2$, $c=2.5$, $d=0.2$, $k=1$ and $K=1$ which satisfy the bifurcation criteria of Theorem 2. The equilibrium point is then $P_1(0.92,0.75,-0.37,-0.14)$. The initial conditions are taken as $x(0)=y(0)=z(0)=u(0)=2$. The critical delay value is $\tau=0.30329$, calculated from \eqref{taucr2}. Again, numerical results given in (a), (b) and (c) are examples for  stable, Hopf bifurcation and unstable cases, respectively. It is also observed that the transition from stability to instability through the Hopf bifurcation is experienced in both $x$ and $y$ variables for the equilibrium point $P_1$. The  results are in well accordance  with predictions of our theoretical results.

After these numerical simulations, we present two and three dimensional phase portraits of the system (Figures 3-6) corresponding to the cases summarized above and they are placed at the end of the manuscript.

\begin{figure}[!hp]
\centering
\subfigure[]{\label{main:a}\includegraphics[scale=.37]{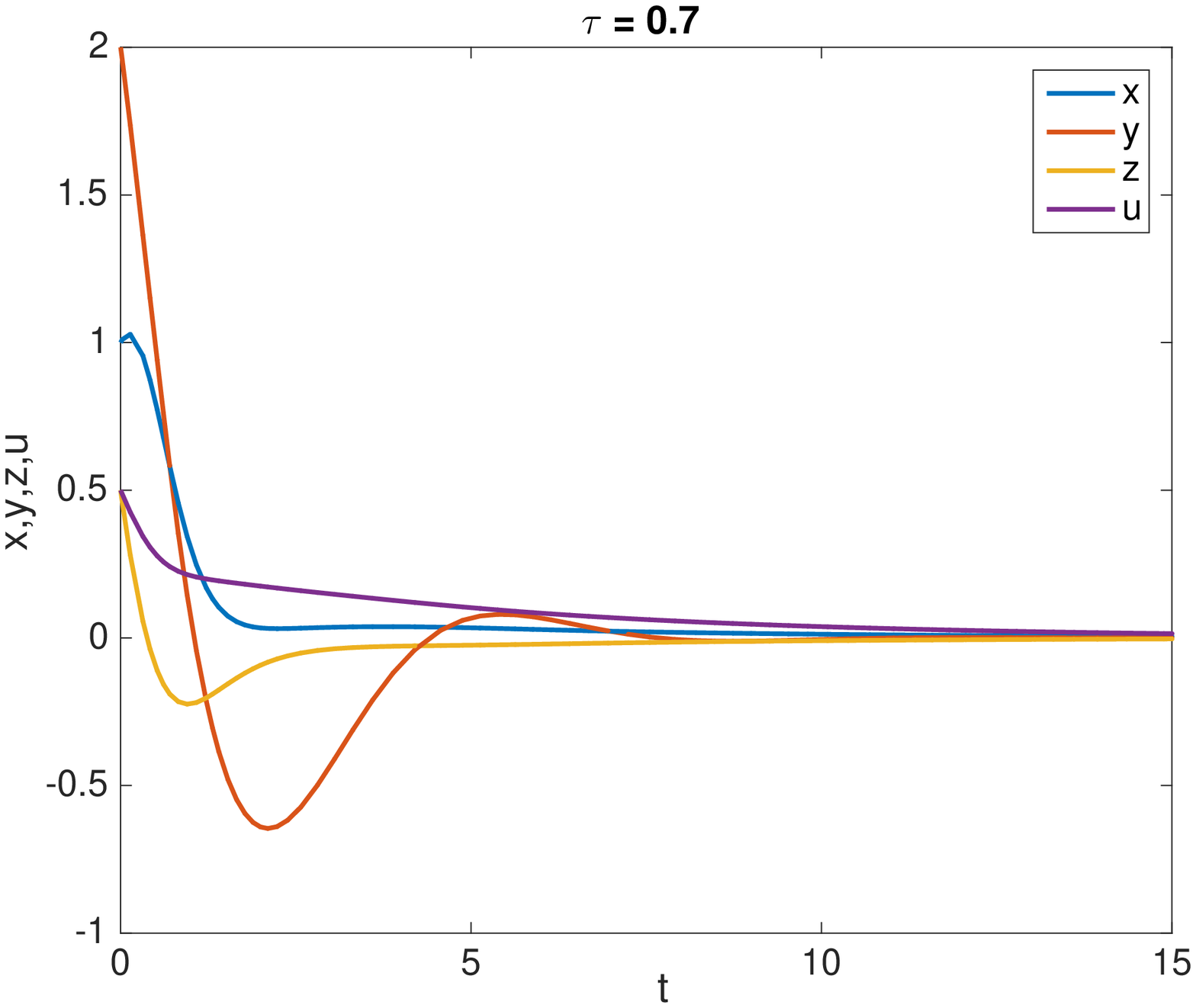}}
\centering
\subfigure[]{\label{main:b}\includegraphics[scale=.50]{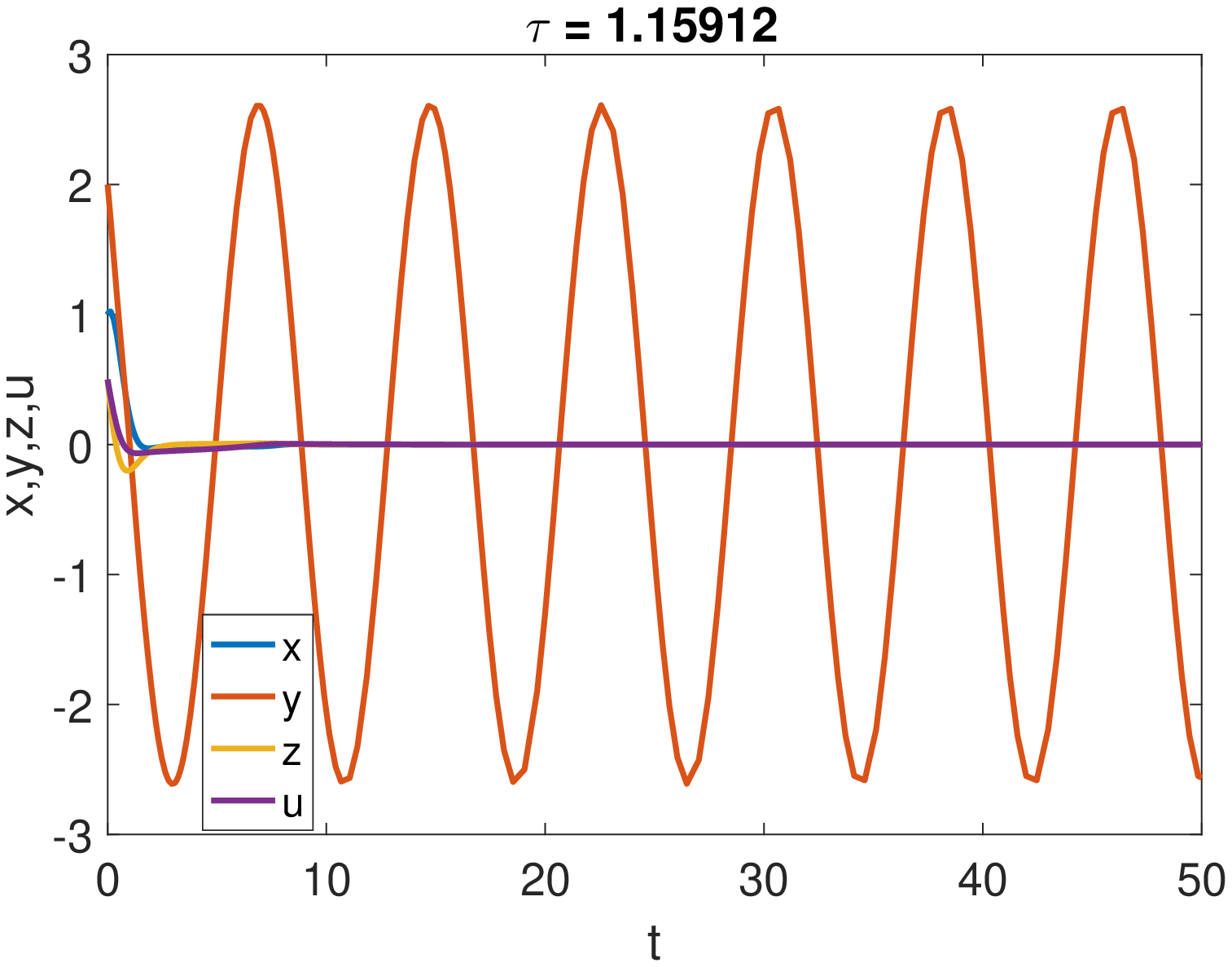}}
\centering
\subfigure[]{\label{main:c}\includegraphics[scale=.45]{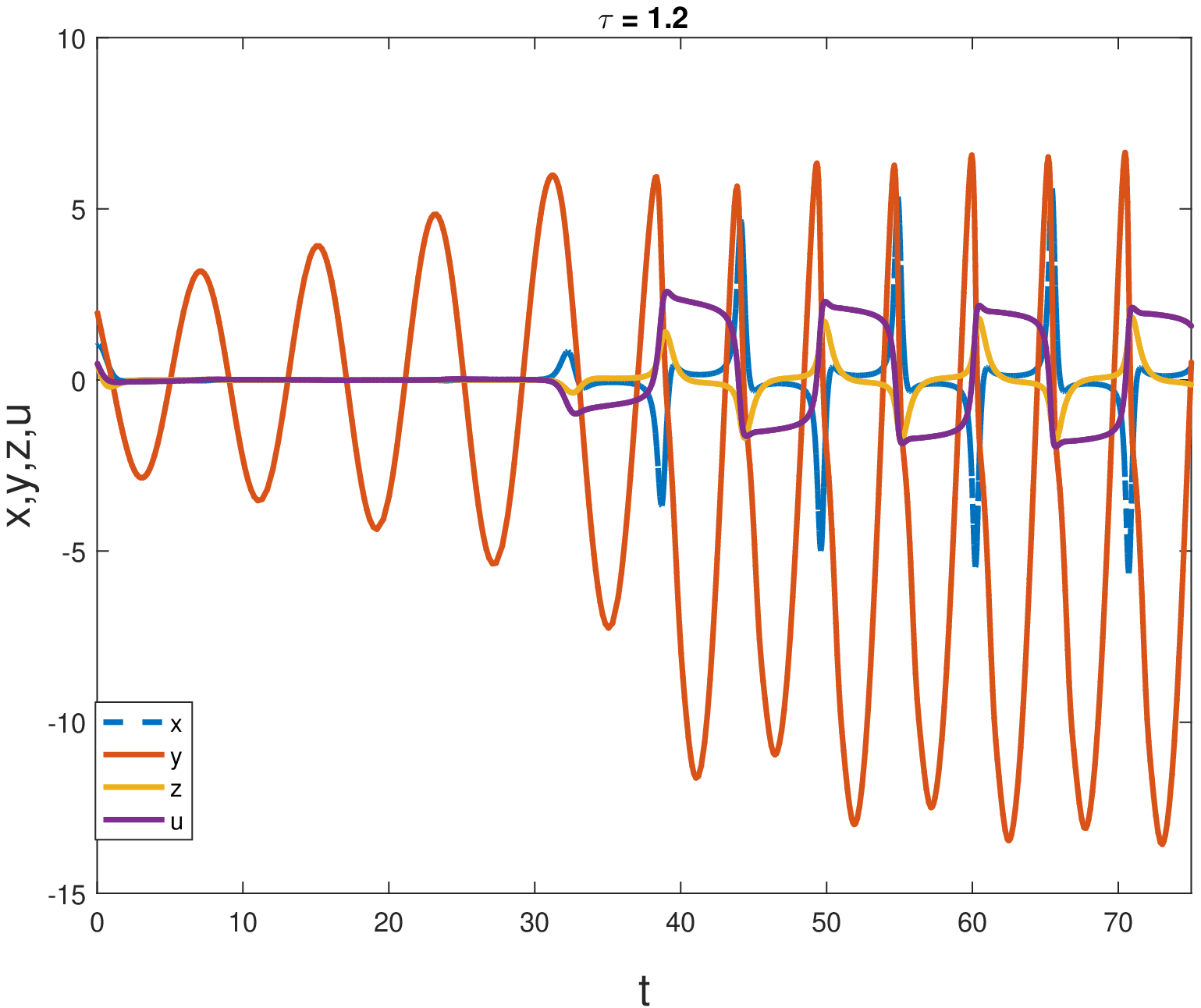}}
\caption{Numerical simulations of the system \eqref{sysX} for different $\tau$ values (a) $\tau=0.7$, (b) $\tau=1.15912$, (c) $\tau=1.2$, for the equilibrium point $P_0$. Here we take the parameter values as $a=5$, $b=0.4$, $c=1.5$, $d=0.2$, $k=0.17$ and $K=1$. Also, the initial conditions are taken as $X(0)=1$, $Y(0)=2$, $Z(0)=U(0)=0.5$. The transition from stability to instability through the Hopf bifurcation is experienced only in the $y$ variable for the equilibrium point $P_0$. }
\label{fig:main}
\end{figure}

\begin{figure}[!hp]
\centering
\subfigure[]{\label{main:a}\includegraphics[scale=.37]{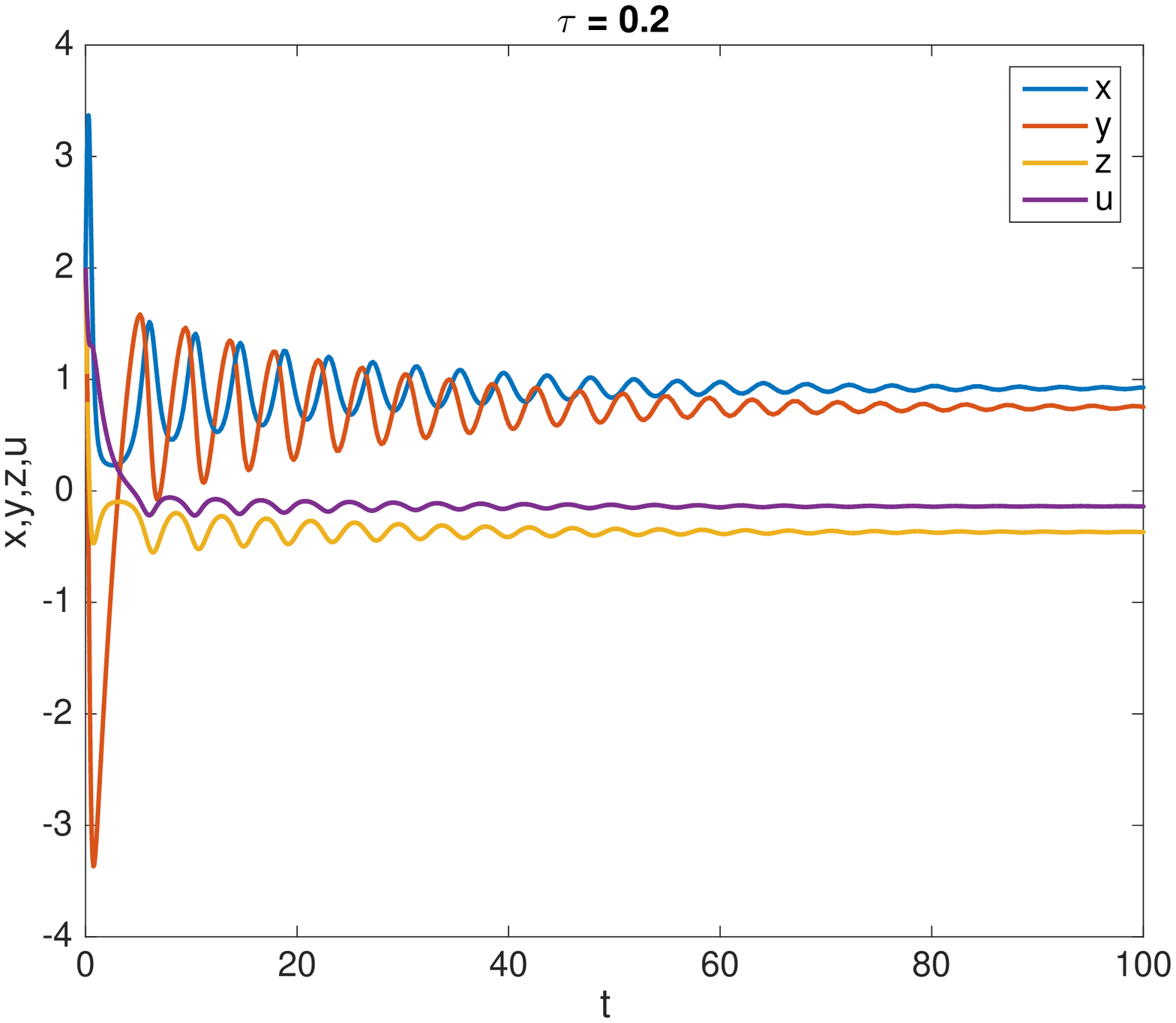}}
\centering
\subfigure[]{\label{main:b}\includegraphics[scale=.38]{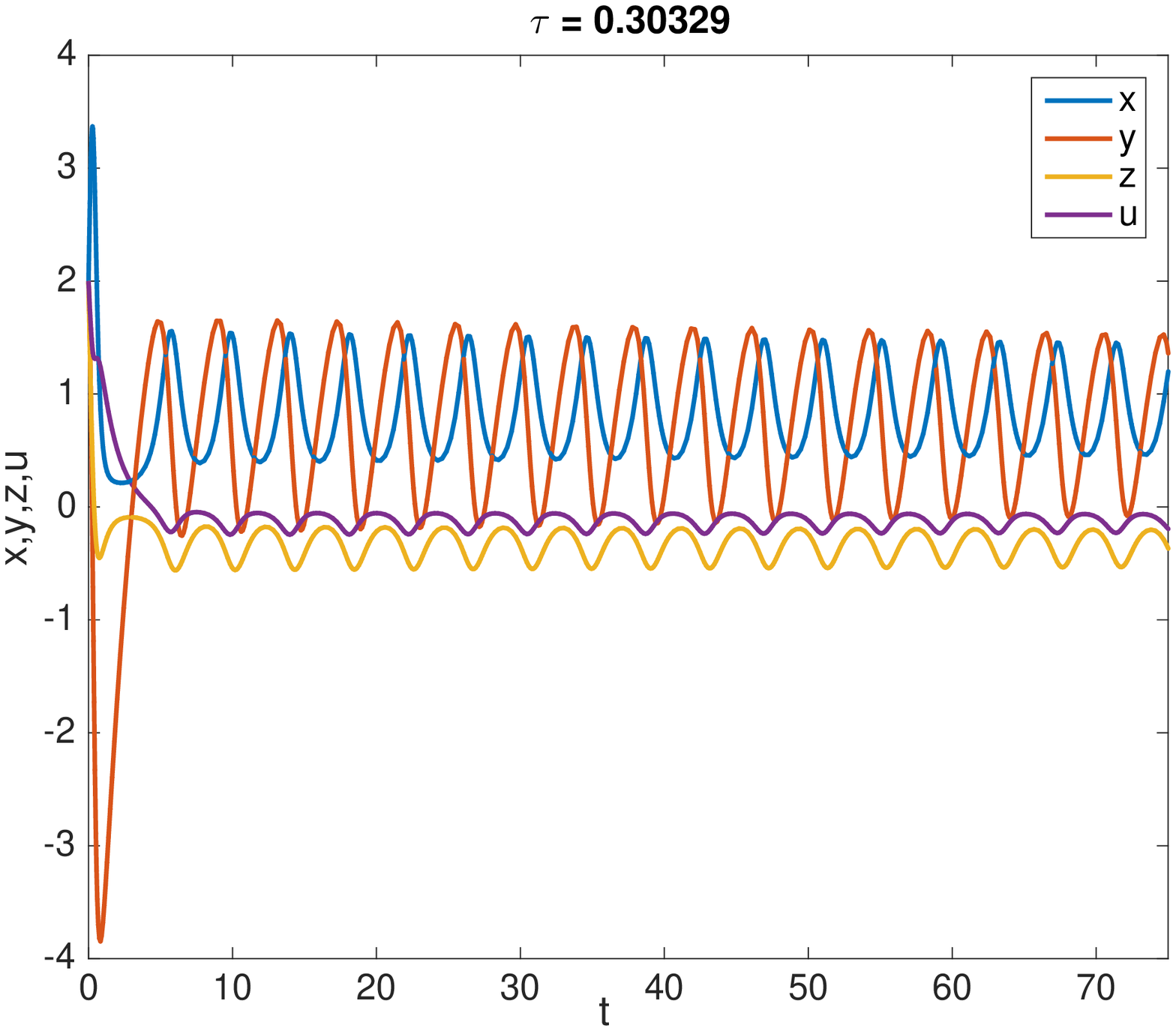}}
\centering
\subfigure[]{\label{main:c}\includegraphics[scale=.37]{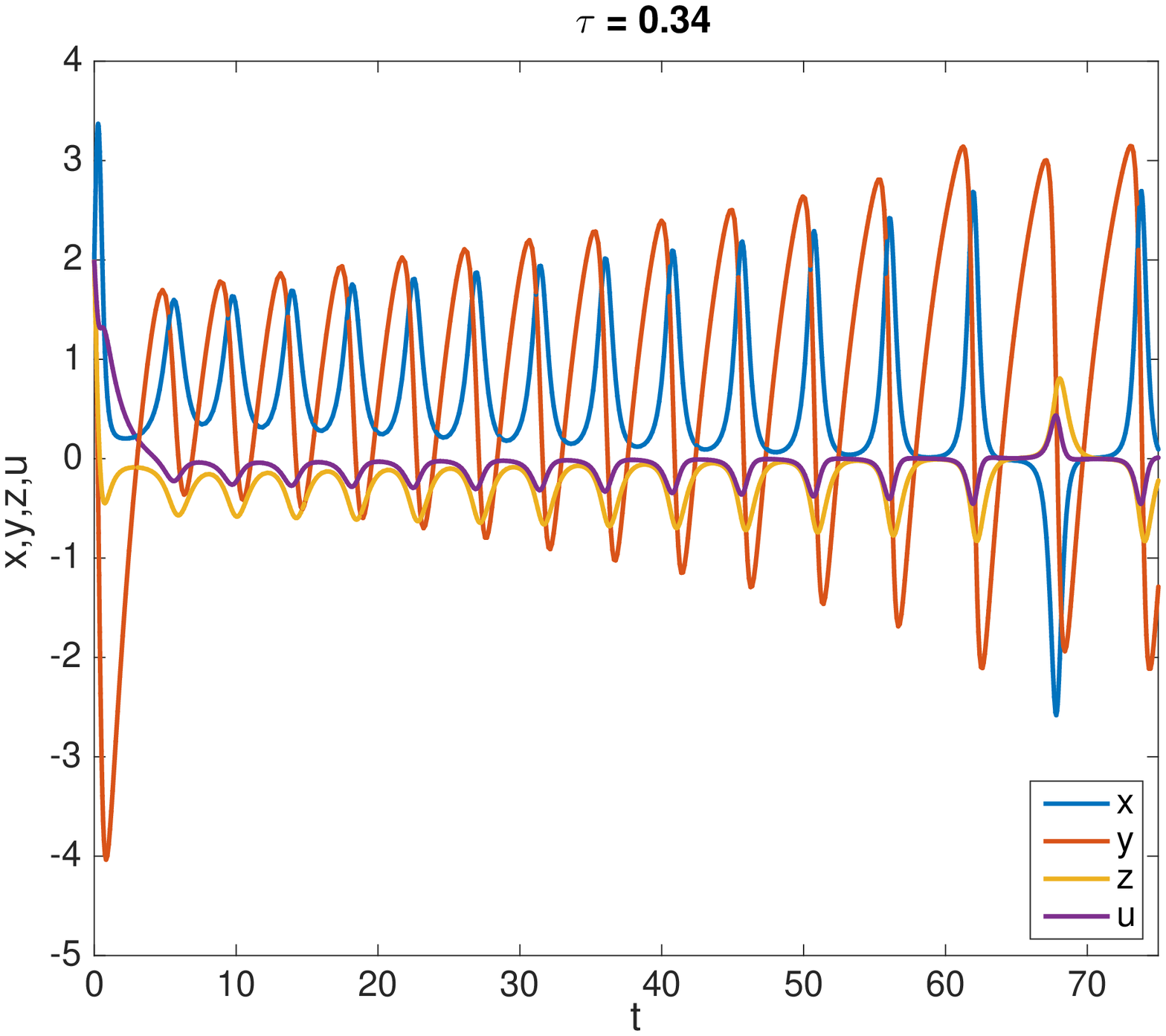}}
\caption{Numerical simulations of the system \eqref{main} for different $\tau$ values (a) $\tau=0.2$, (b) $\tau=0.30329$, (c) $\tau=0.34$, for the critical point $P_1$. Here we take the parameter values as $a=0.2$, $b=0.2$, $c=2.5$, $d=0.2$, $k=1$ and $K=1$. Also, the initial conditions are taken as $x(0)=y(0)=z(0)=u(0)=2$. The transition from stability to instability through the Hopf bifurcation is experienced in both $x$ and $y$ variables for the equilibrium point $P_1$. }
\label{fig:main}
\end{figure}

\section{Conclusion}
In this work, a new dynamical finance system is established. The new system's basic dynamical behaviors,  stability and Hopf bifurcation are investigated at the equilibrium points. We analysed the system \eqref{main}, which we constructed on two existing models in the literature. In the system $S_1$, which is given in \eqref{main1}, there are three state variables, $x,y,z$ and $S_1$ includes a delay term in the variable $y$. In the system $S_2$, which is in \eqref{main2}, there are four state variables, $x,y,z,u$. When $K=0$ in $S_1$ and $u=0$ in $S_2$, they coincide and become the system \eqref{eq60}. Our main system $S_m$, Eq. \eqref{main} is a composition of $S_1$ and $S_2$, with four state variables $x,y,z,u$, the delay feedback coefficient $K$, and the parameters $a,b,c,d,k$. Since $S_m$ is obtained by adding a delay term to $S_2$, it reflects the delay effect on the system $S_2$.

Time delay parameter $\tau$ is taken as a bifurcation and control parameter in order to search the system's stability behaviour. After linearization, the characteristic equations are examined at the equilibrium points and we proved that a Hopf bifurcation exists. If time delay $\tau$ passes a critical value, the system experiences a Hopf bifurcation,  the stability condition of the system changing from stable to unstable. Through numerical simulations, our main results  are confirmed; that the system undergoes a Hopf bifurcation with appropriate parameters and some graphs are shown at different time delay  values.

The equilibrium points of $S_2$ and $S_m$ are the same. $S_m$ differs from the system  $S_2$ by the delay feedback term. For the values of the parameters considered in Section 3.2, $S_m$ when $K=0$, namely the system $S_2$  is stable. When $K=1$, $S_m$ is stable for some range of the time-delay term $\tau$; however, it becomes unstable after this critical threshold. Therefore, we can say that, there are cases in which time-delay term has a destabilizing effect on $S_2$.

We worked out bifurcation analysis of a  dynamic finance system and we found that the system has rich dynamic behaviors and responses. Then, this study can be helpful for the relevant fields, especially for the financial models, as a theoretical reference and it  deserves to be studied more. As an open problem for further investigation, for instance, we can mention the search for   chaotic or hyperchaotic character of the system.

\bibliographystyle{spmpsci}

\newpage

\begin{figure}[!hp]
\subfigure[]{\label{sek21} 
\includegraphics[scale=0.41]{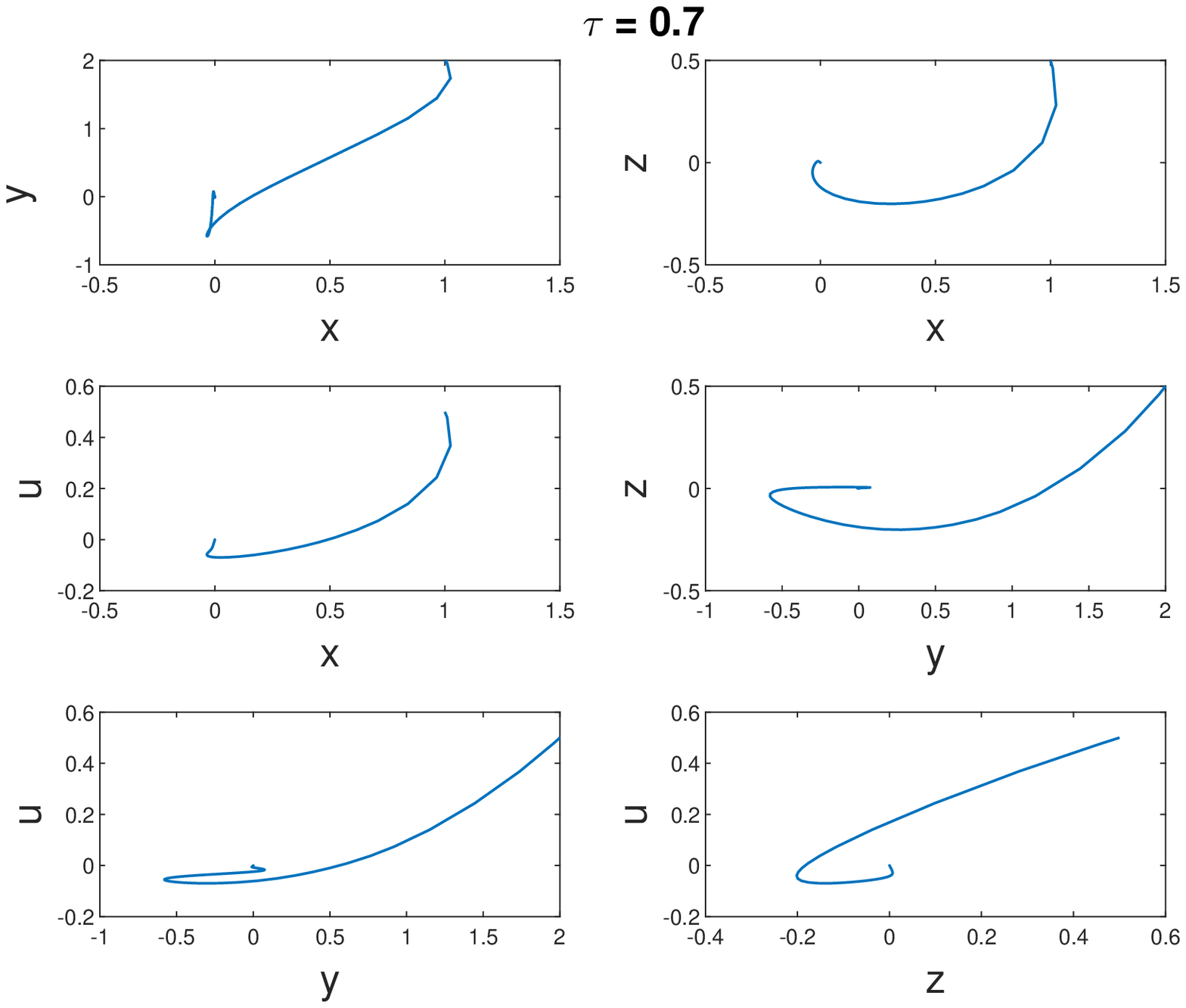}}
\subfigure[]{\label{sek23} 
\includegraphics[scale=0.35]{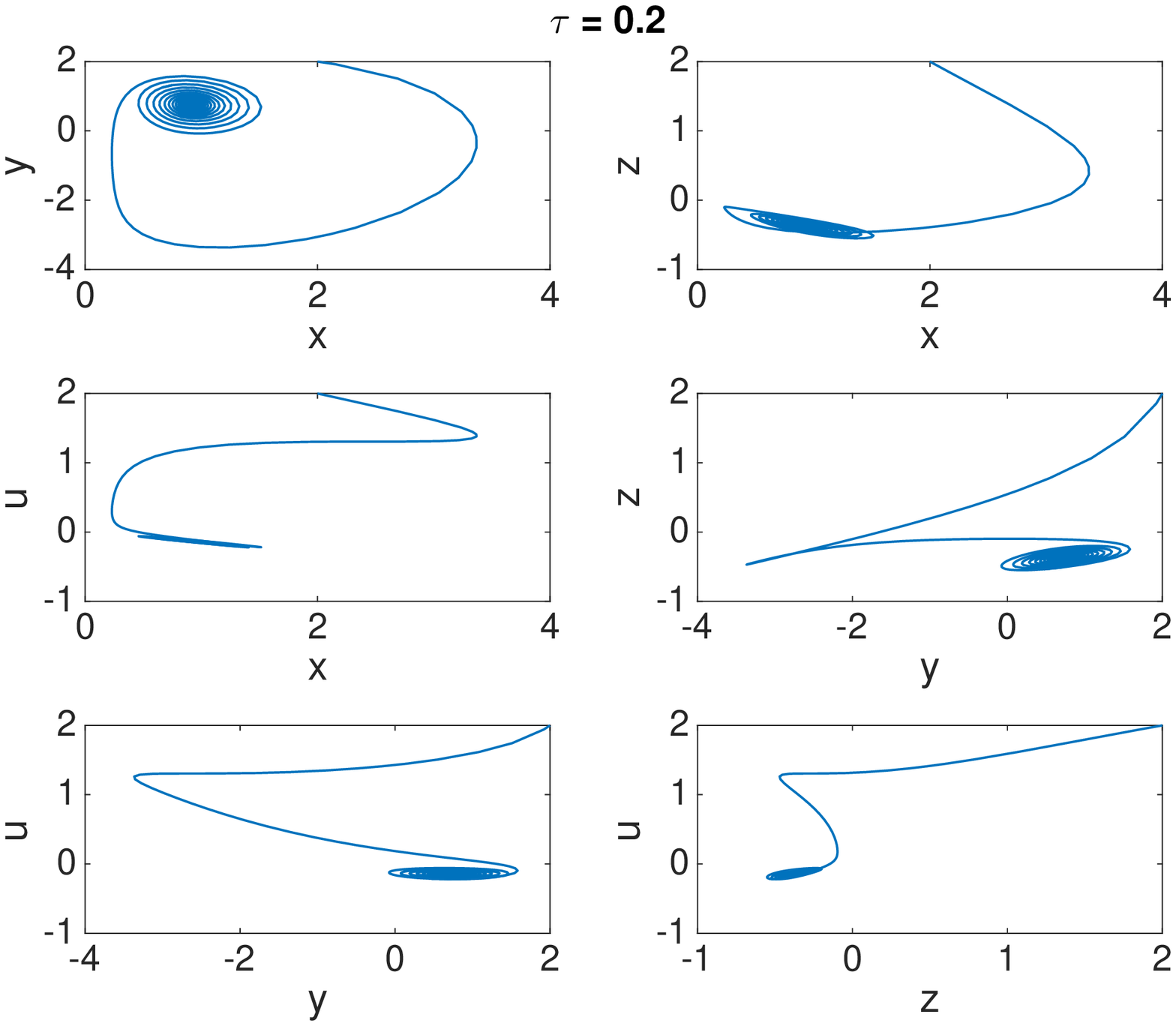}}
\subfigure[]{\label{sek22} 
\includegraphics[scale=0.41]{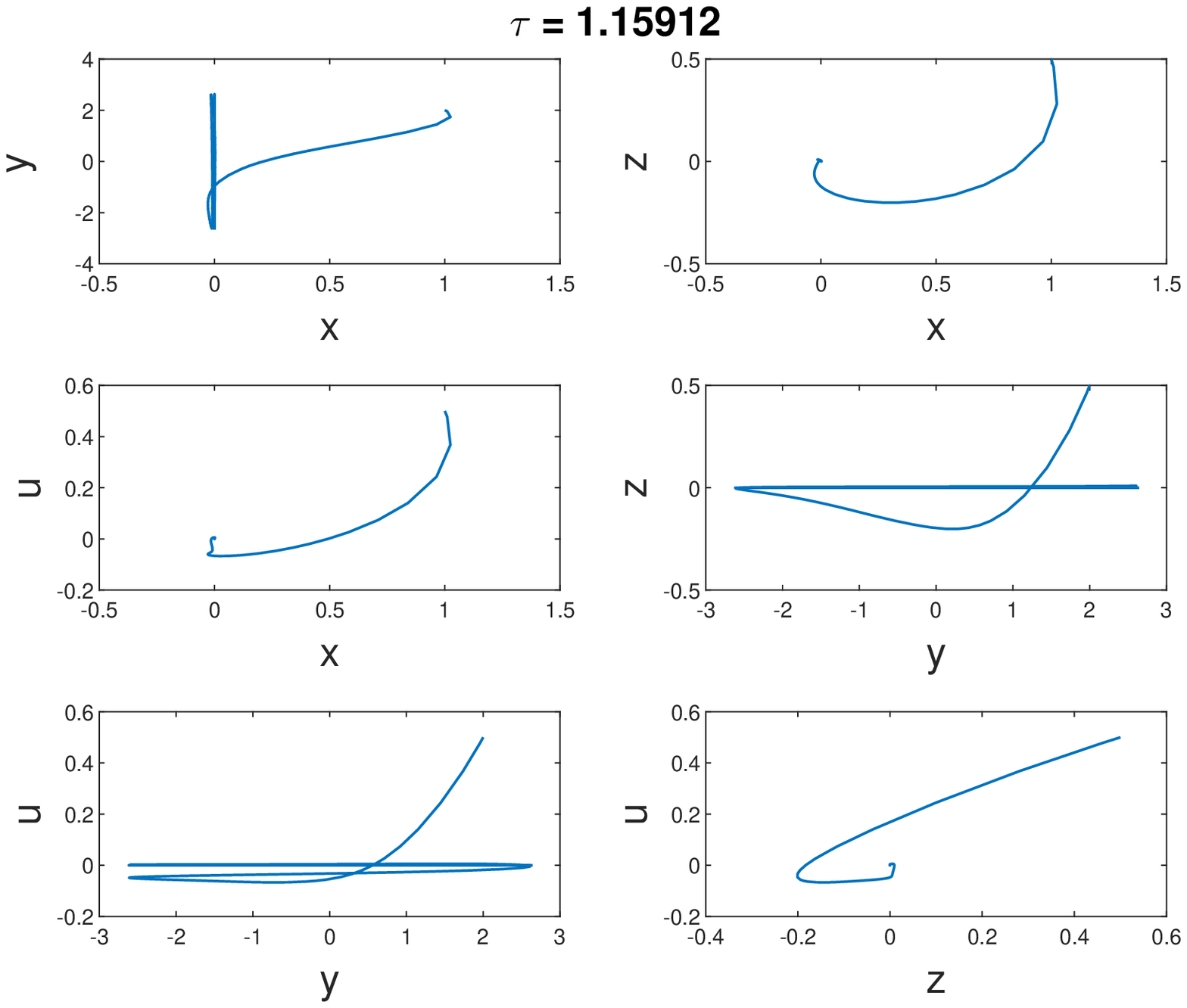}}
\subfigure[]{\label{sek23} 
\includegraphics[scale=0.35]{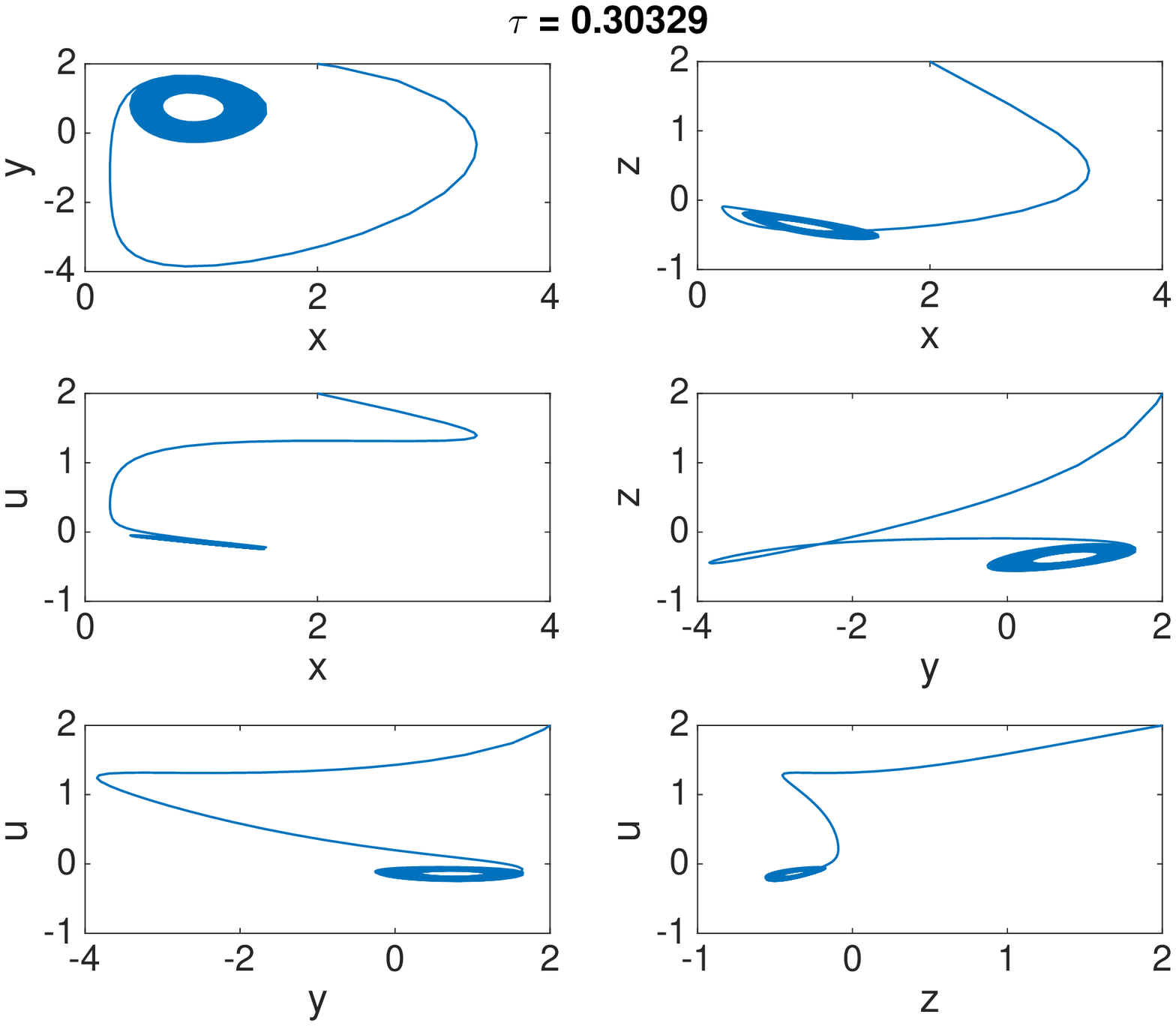}}
\subfigure[]{\label{sek21} 
\includegraphics[scale=0.39]{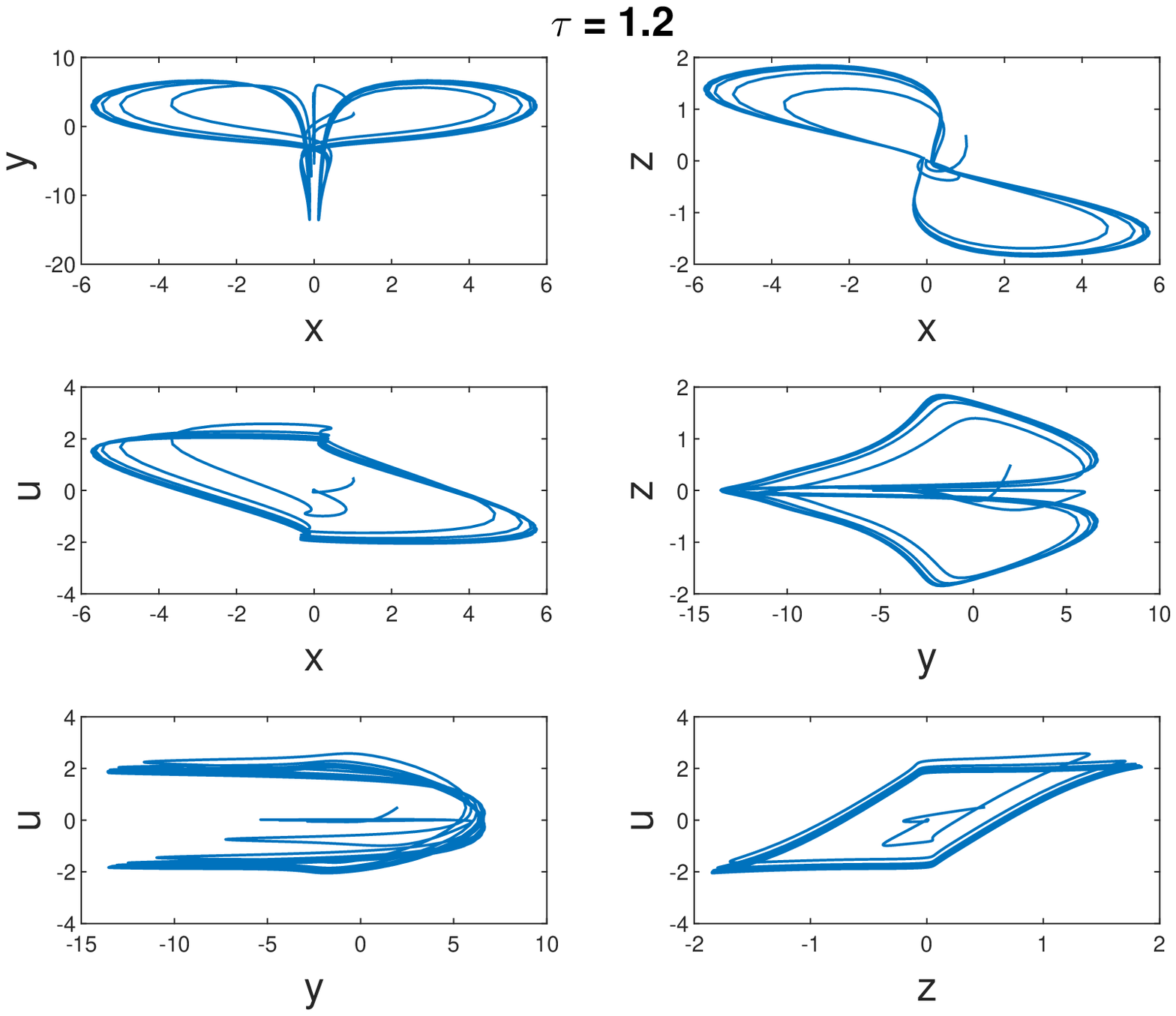}}
\subfigure[]{\label{sek22} 
\hspace{1cm}
\includegraphics[scale=0.35]{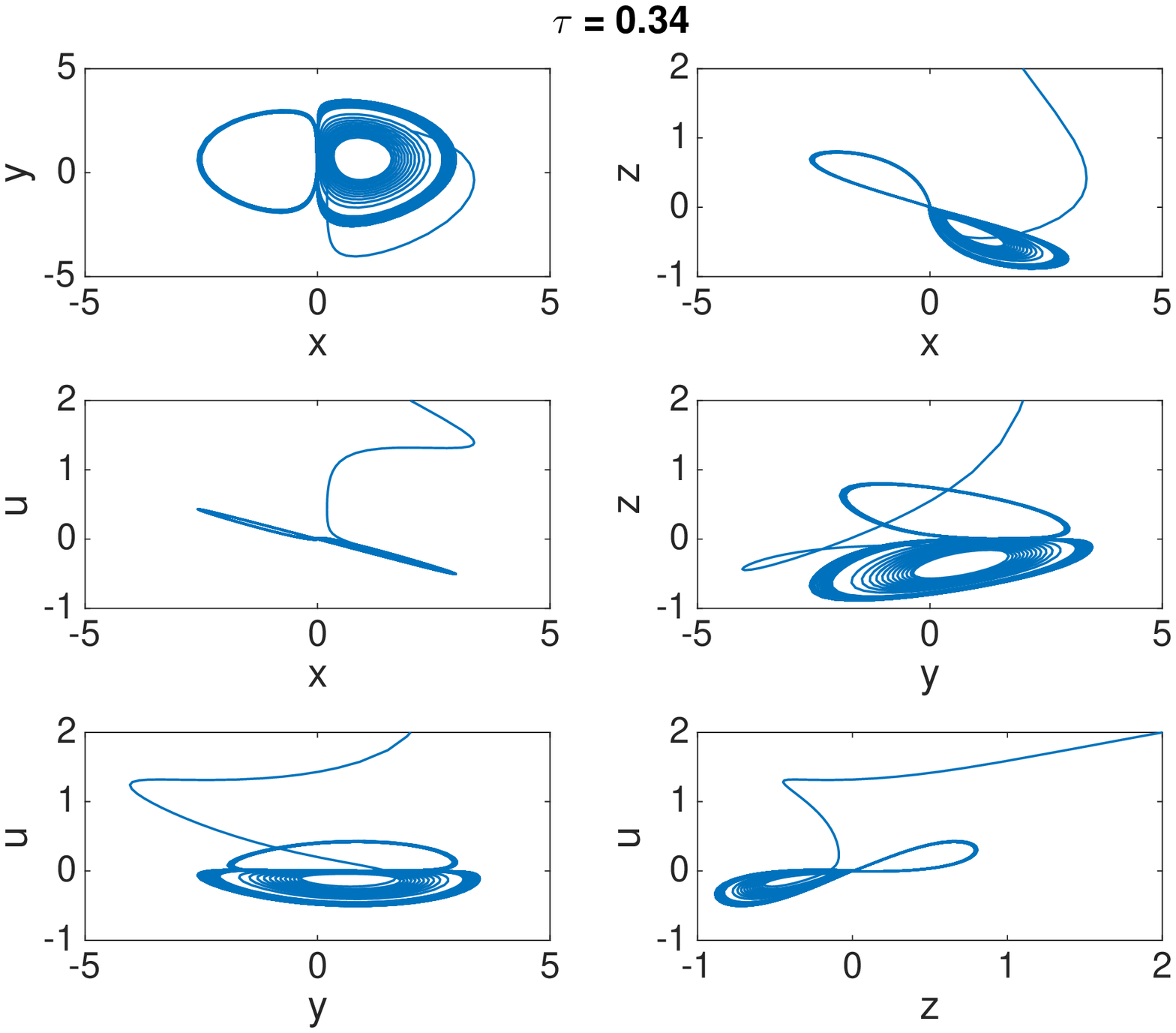}}
\caption{Two dimensional phase portraits obtained from the numerical solutions of the system \eqref{main} for the chosen parameter values and initial conditions; (a), (c), (e) for $P_0$, (b), (d), (f) for $P_1$.}
\label{imper} 
\end{figure}

\begin{figure}[!hp]
\subfigure[]{\label{sek21} 
\includegraphics[scale=0.40]{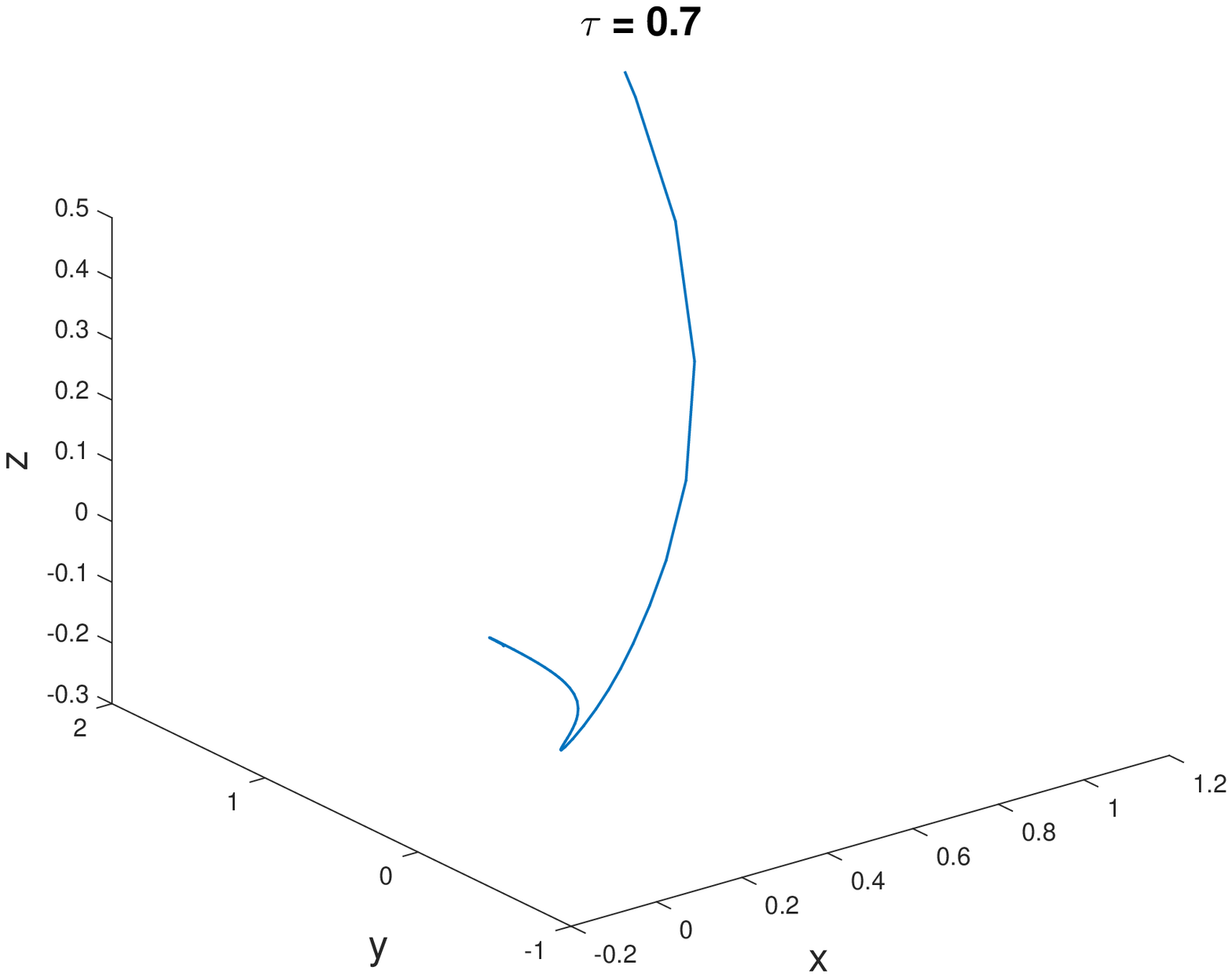}}
\subfigure[]{\label{sek21} 
\includegraphics[scale=0.40]{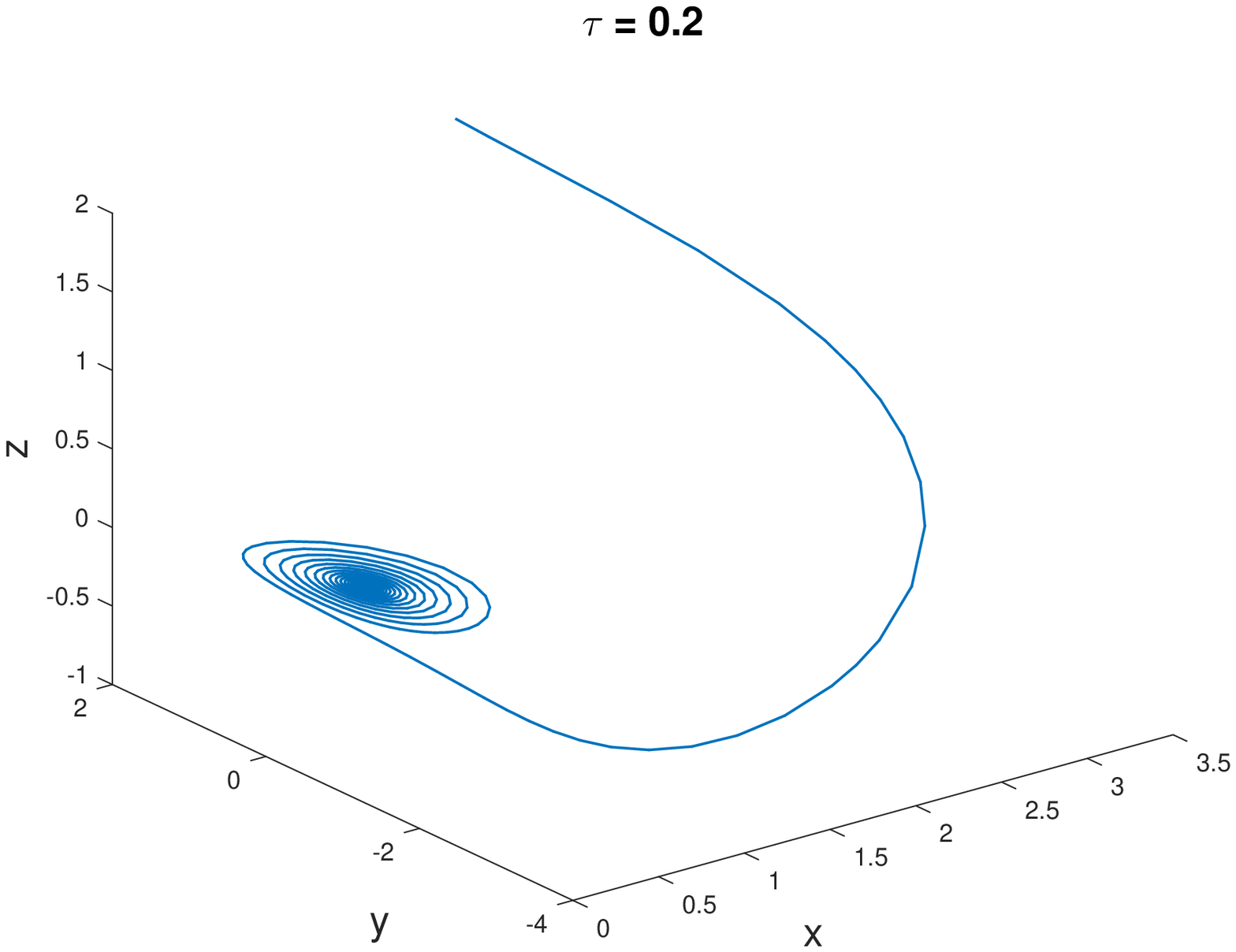}}
\subfigure[]{\label{sek22} 
\includegraphics[scale=0.40]{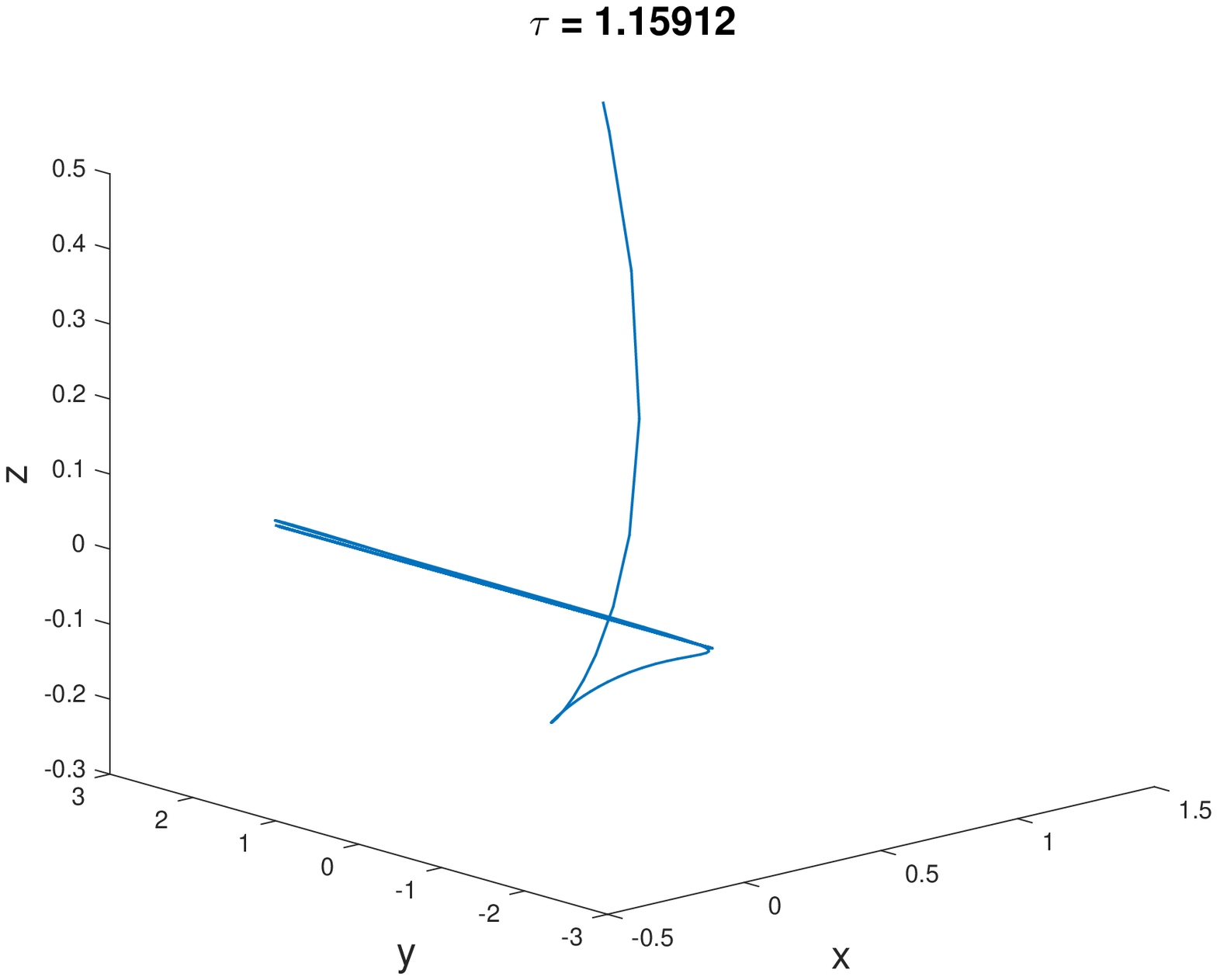}}
\subfigure[]{\label{sek22} 
\includegraphics[scale=0.37]{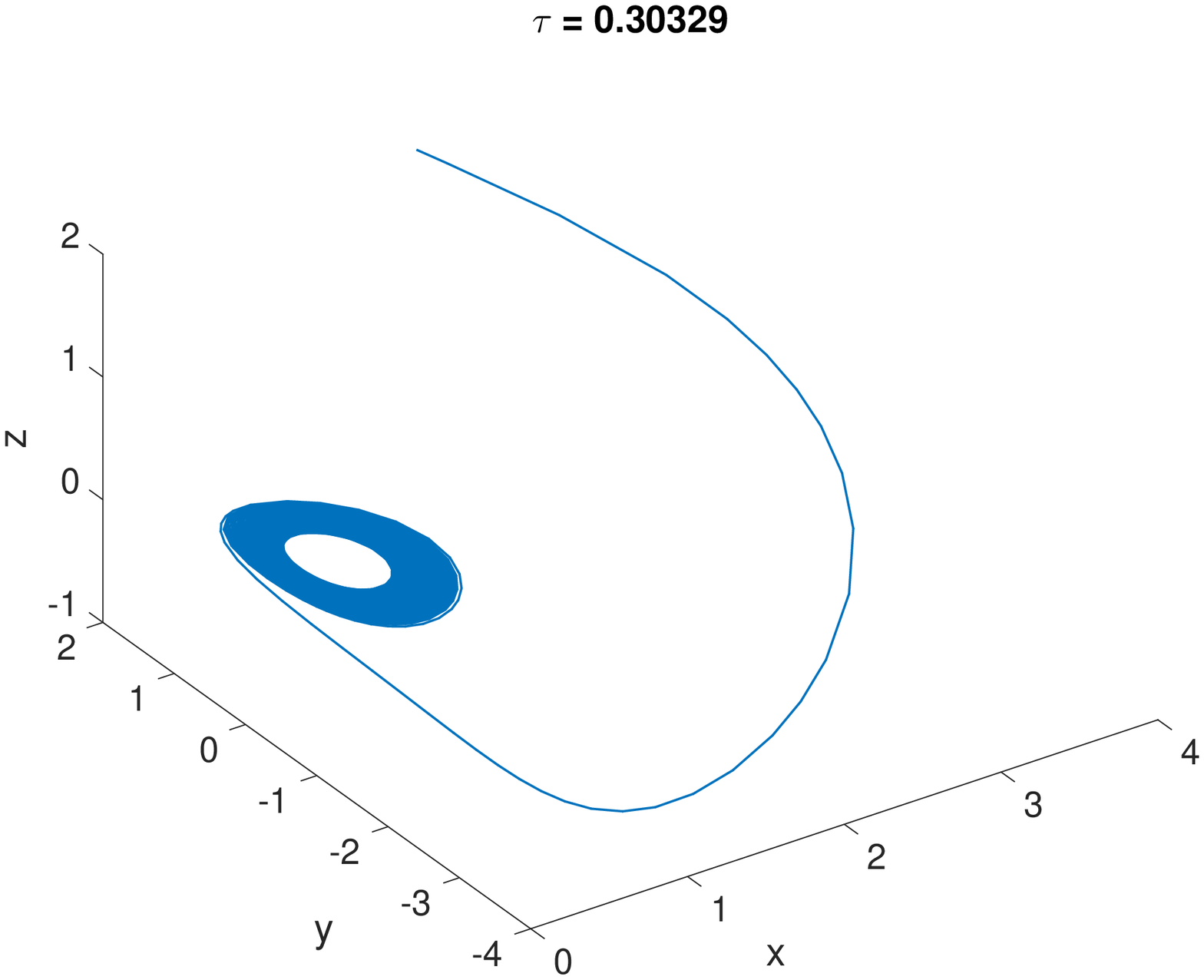}}
\subfigure[]{\label{sek23} 
\includegraphics[scale=0.42]{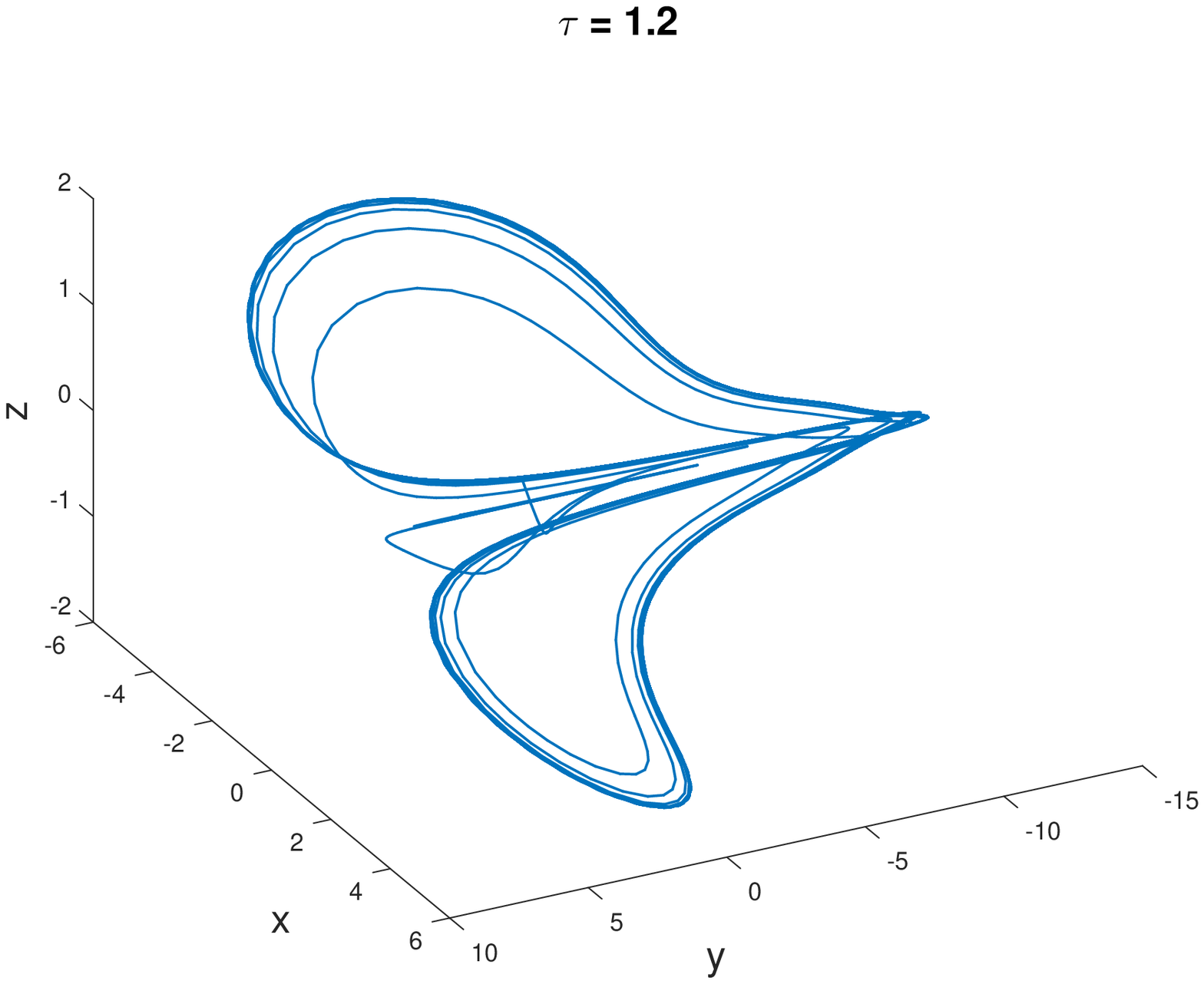}}
\subfigure[]{\label{sek23} 
\hspace{.5cm}
\includegraphics[scale=0.37]{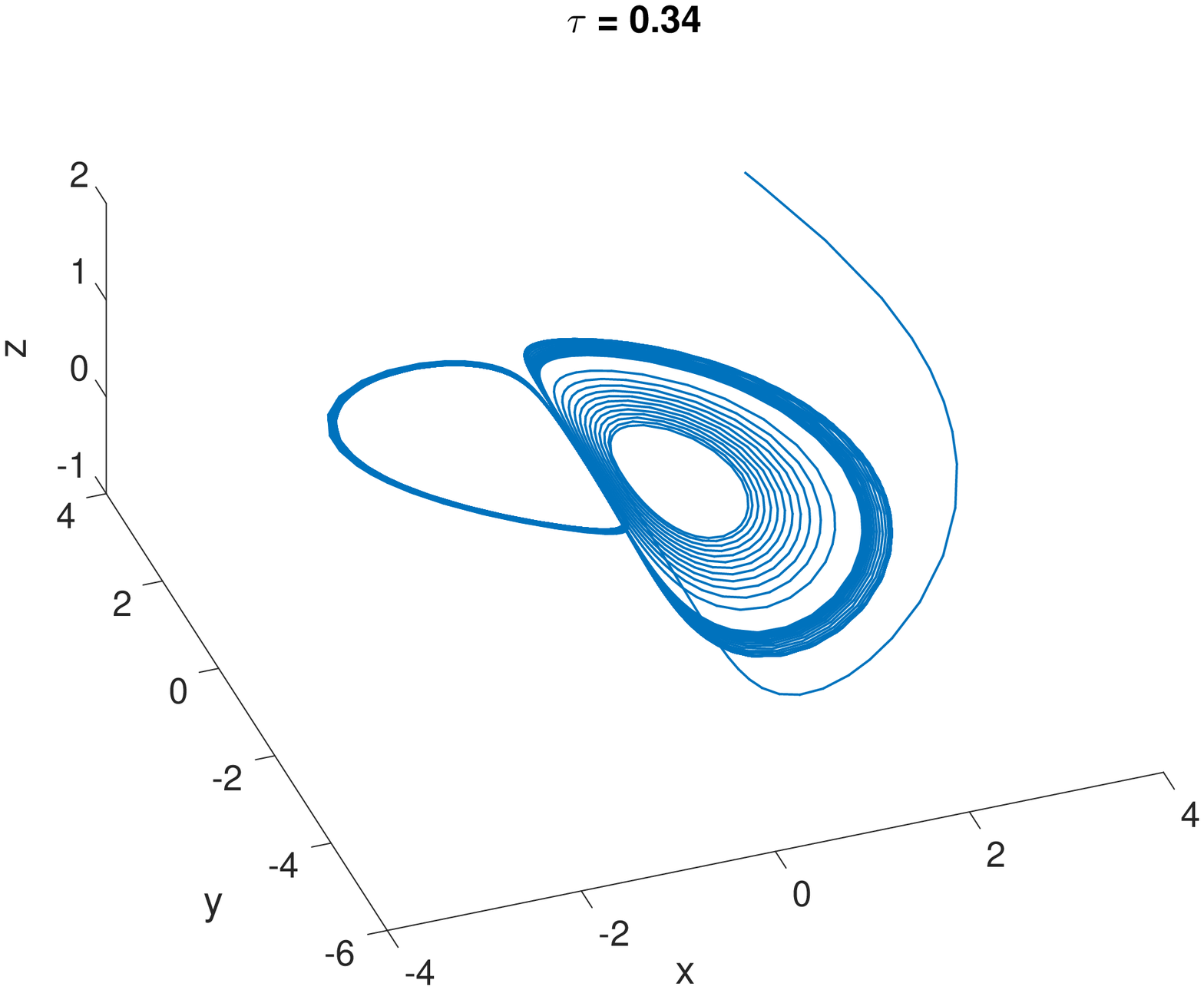}}
\caption{Three dimensional phase portraits of the variables $x$, $y$ and $z$ obtained from the numerical solutions of the system \eqref{main} for the chosen parameter values and initial conditions; (a), (c), (e) for $P_0$, (b), (d), (f) for $P_1$.}
\label{imper} 
\end{figure}

\begin{figure}[!hp]
\subfigure[]{\label{sek21} 
\includegraphics[scale=0.33]{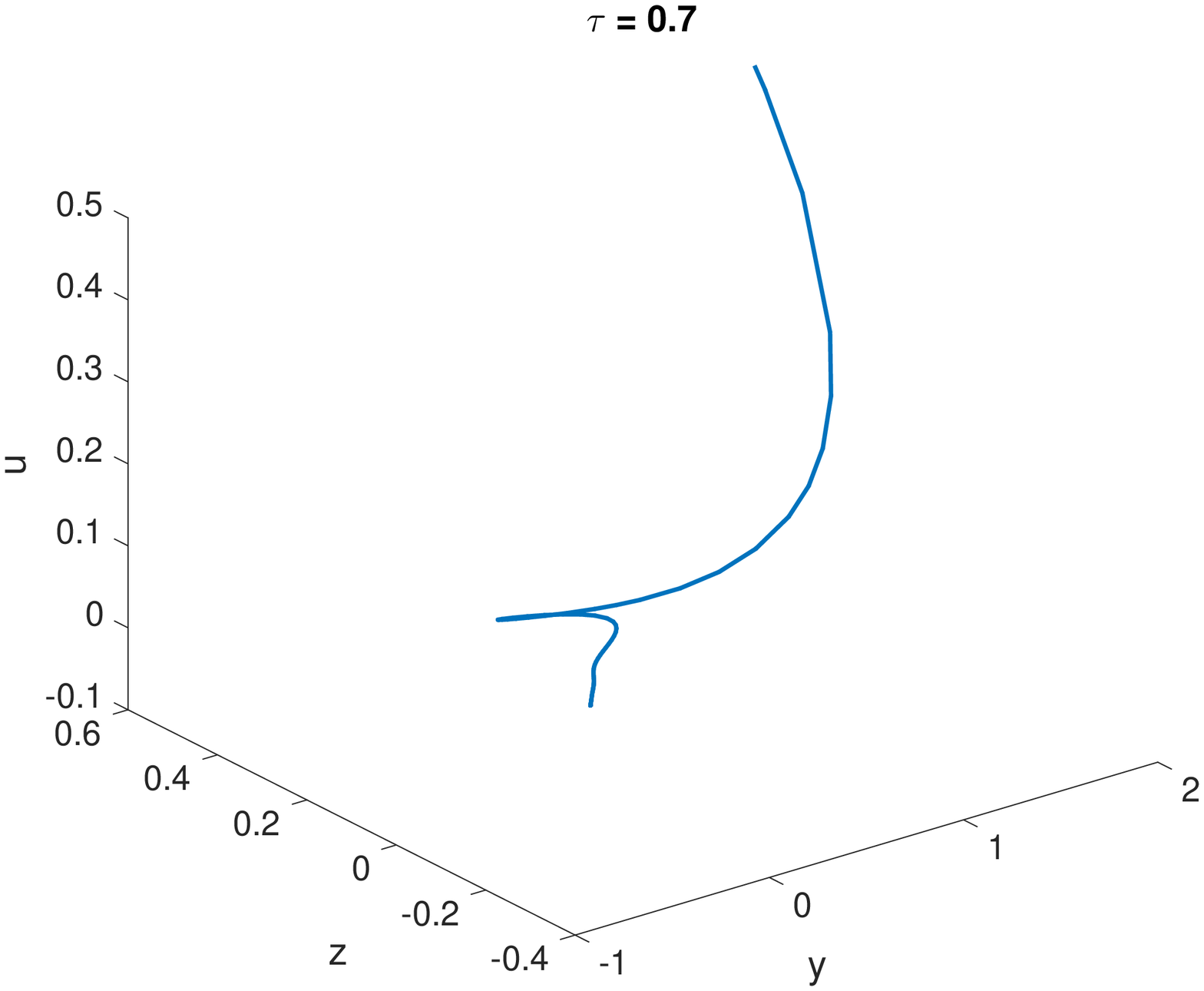}}
\subfigure[]{\label{sek21} 
\includegraphics[scale=0.35]{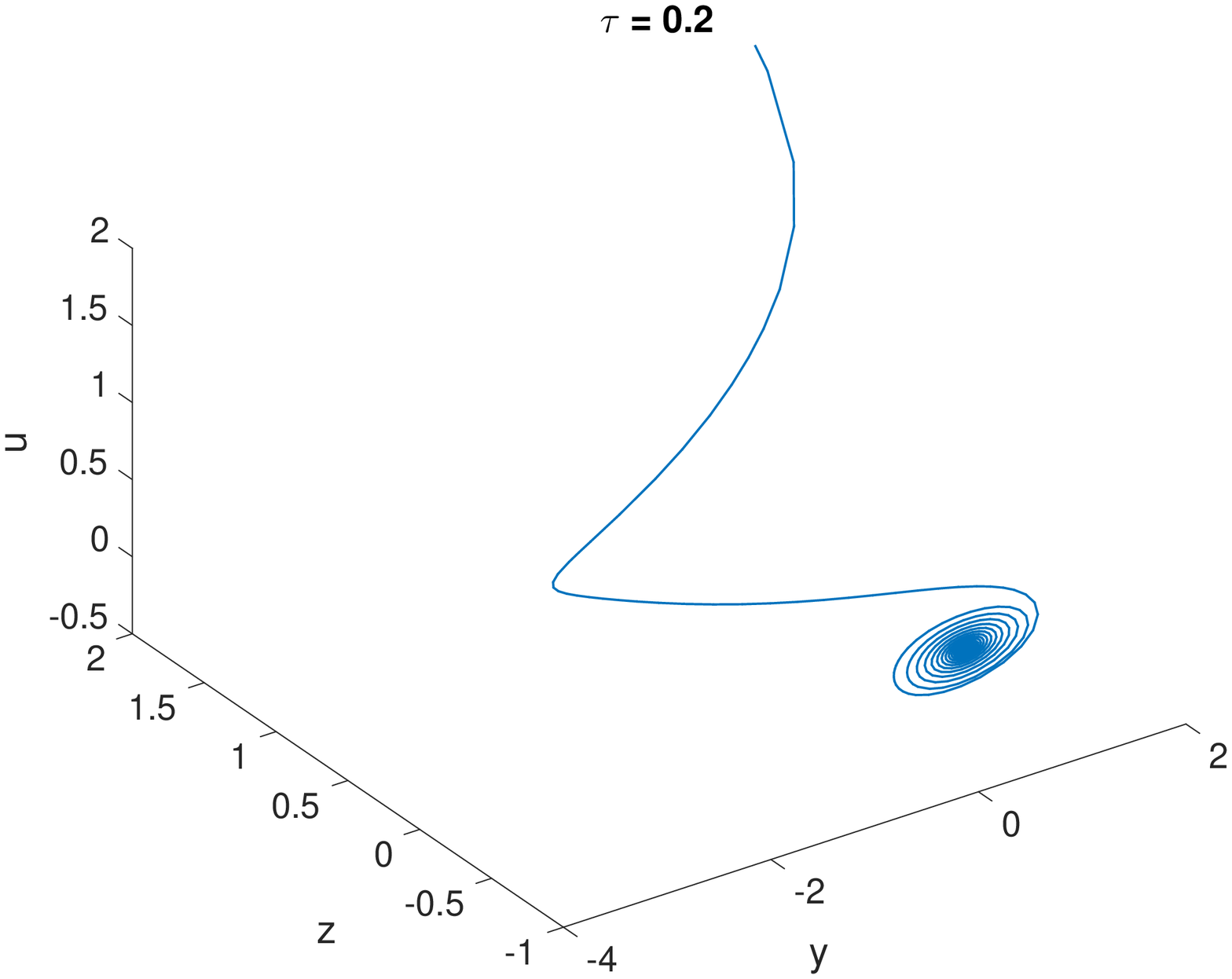}}
\subfigure[]{\label{sek22} 
\includegraphics[scale=0.33]{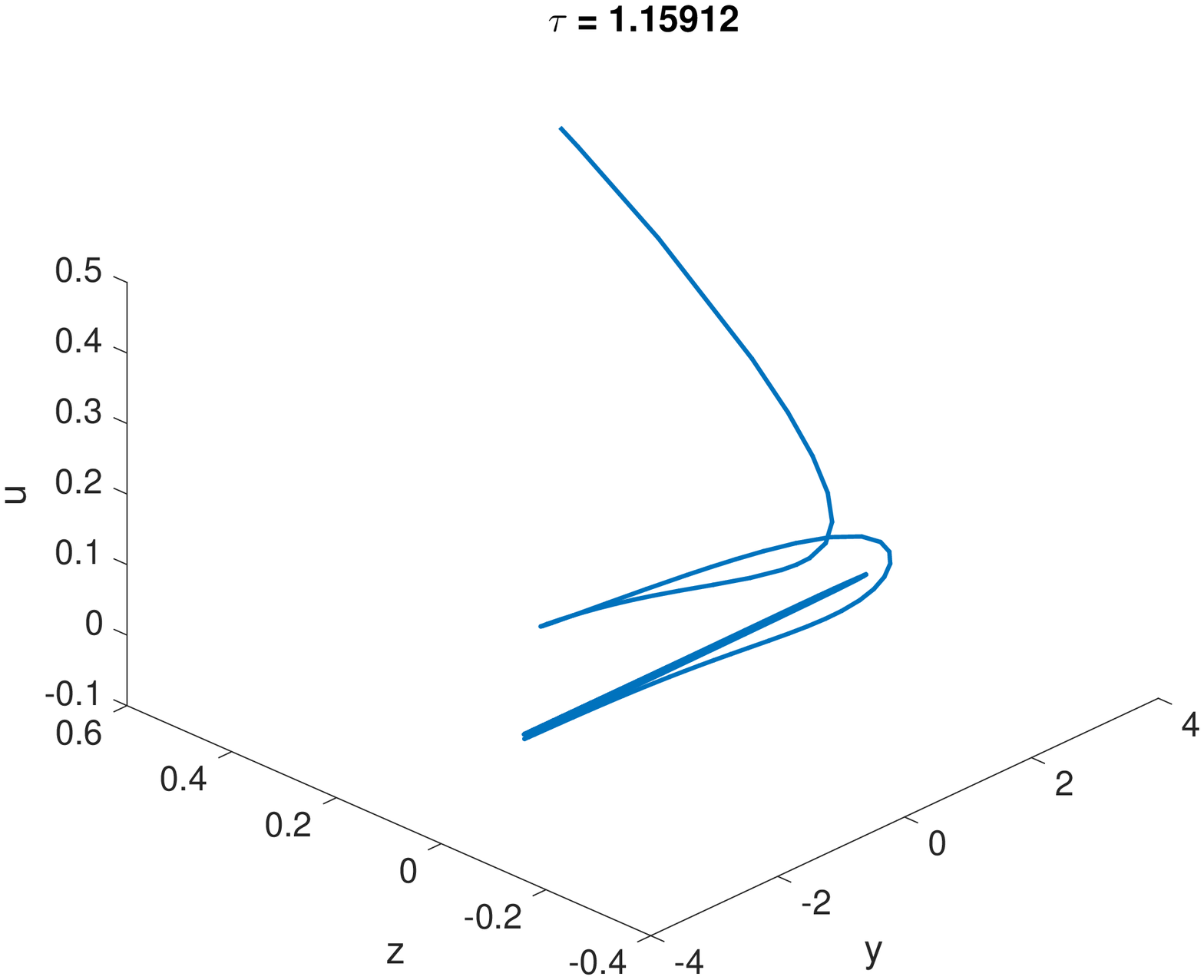}}
\subfigure[]{\label{sek22} 
\includegraphics[scale=0.35]{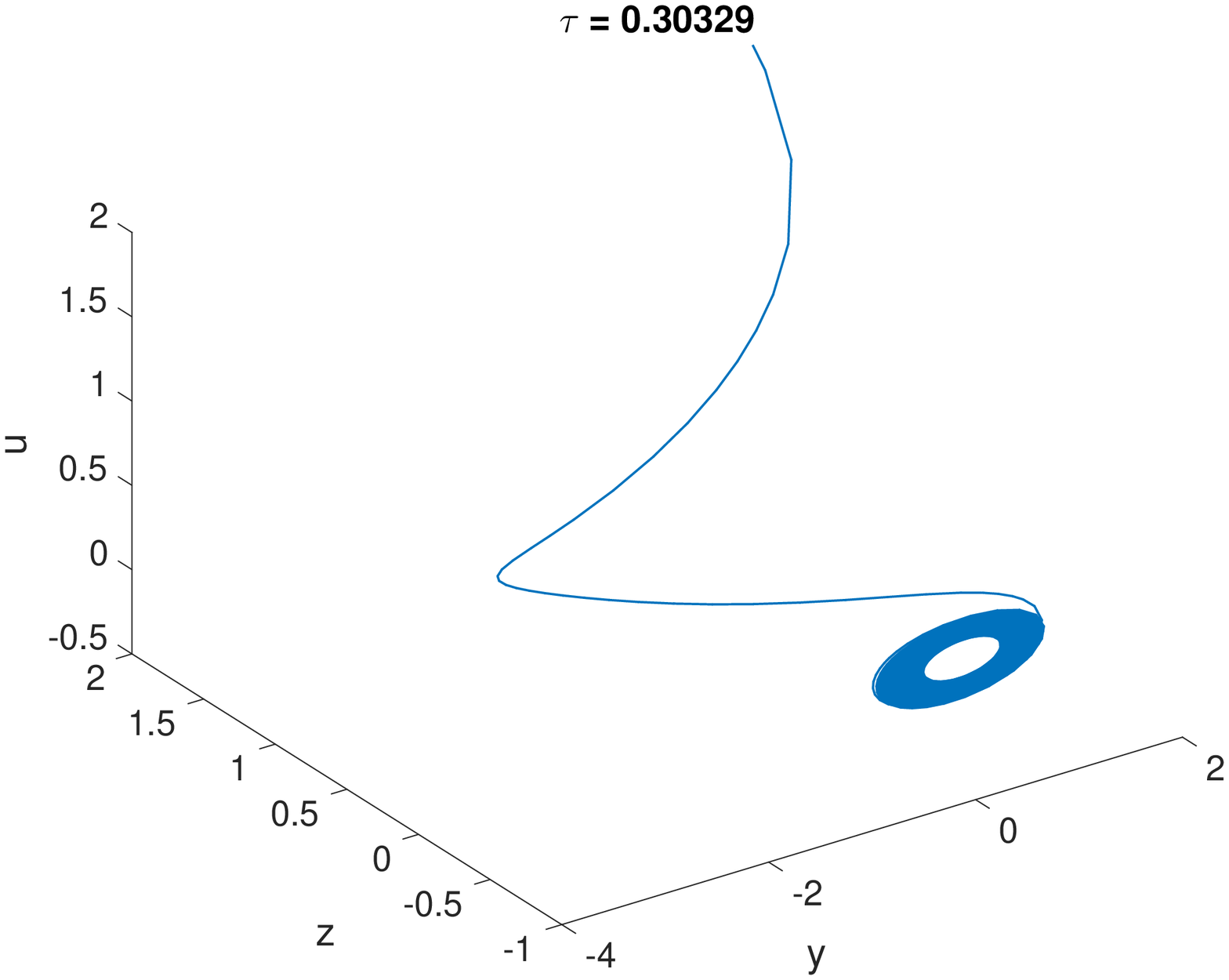}}
\subfigure[]{\label{sek23} 
\includegraphics[scale=0.35]{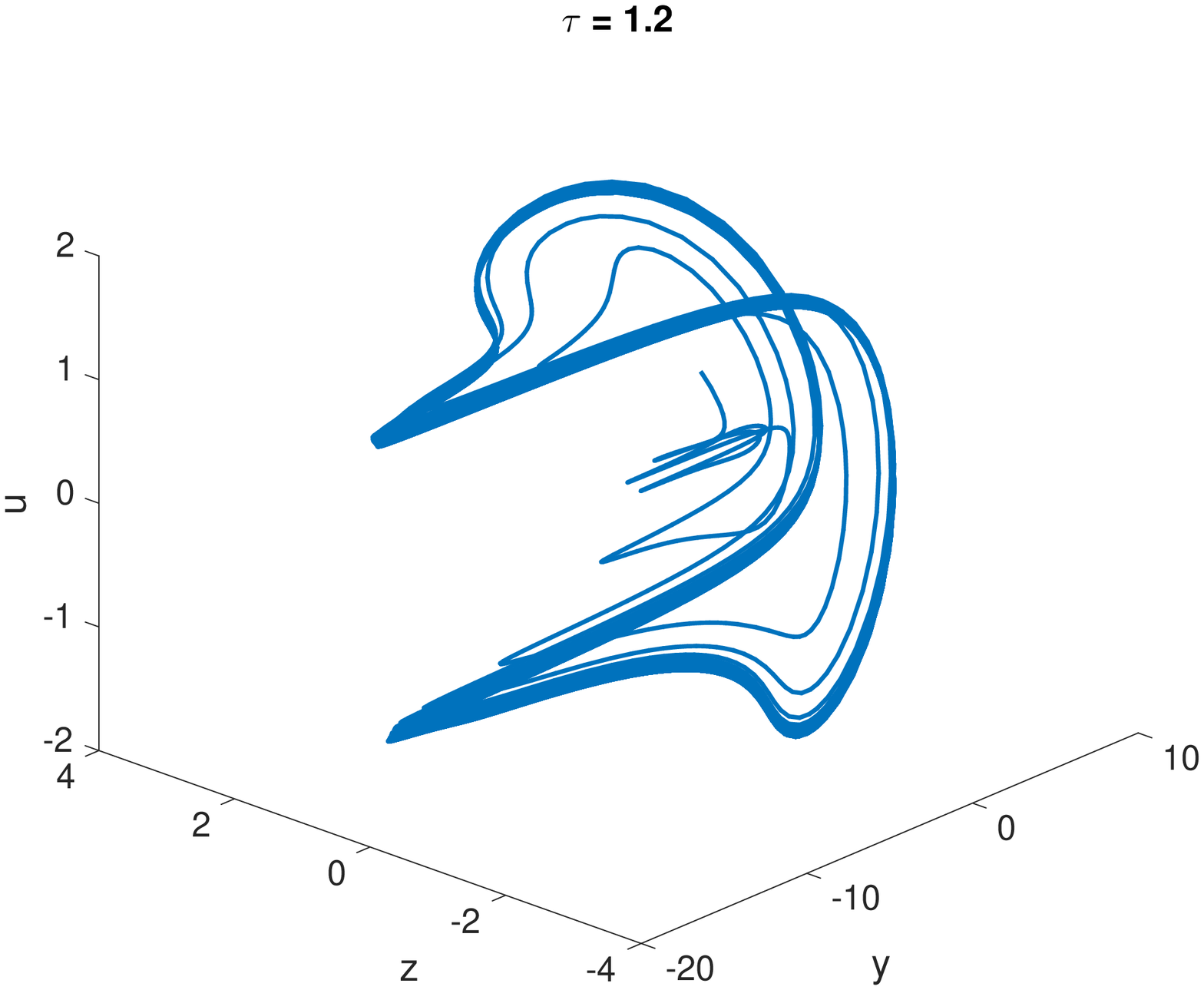}}
\subfigure[]{\label{sek23} 
\hspace{1.5cm}
\includegraphics[scale=0.35]{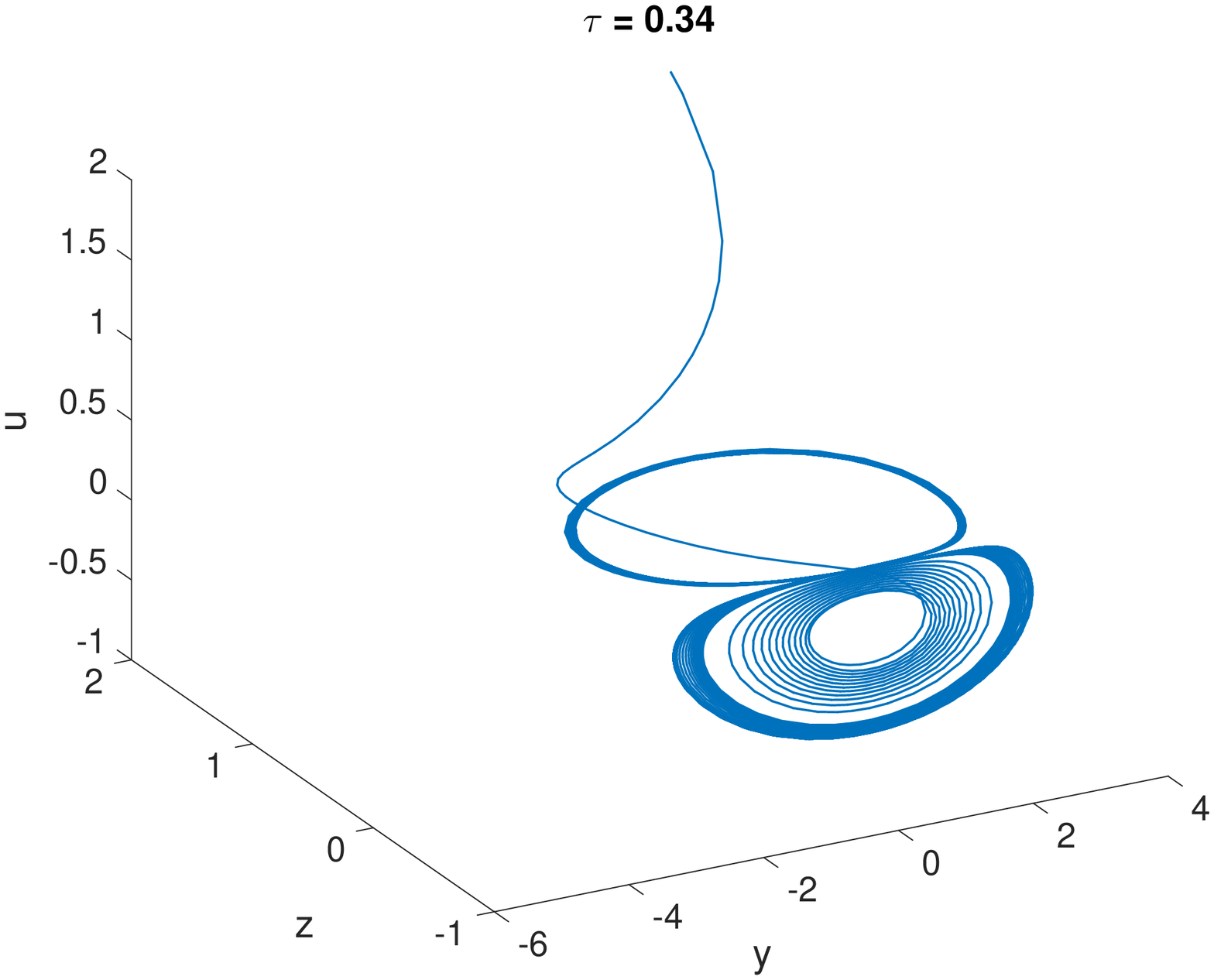}}
\caption{Three dimensional phase portraits of the variables $y$, $z$ and $u$ obtained from the numerical solutions of the system \eqref{main} for the chosen parameter values and initial conditions; (a), (c), (e) for $P_0$, (b), (d), (f) for $P_1$.}
\label{imper} 
\end{figure}



\begin{figure}[!hp]
\subfigure[]{\label{sek21} 
\includegraphics[scale=0.41]{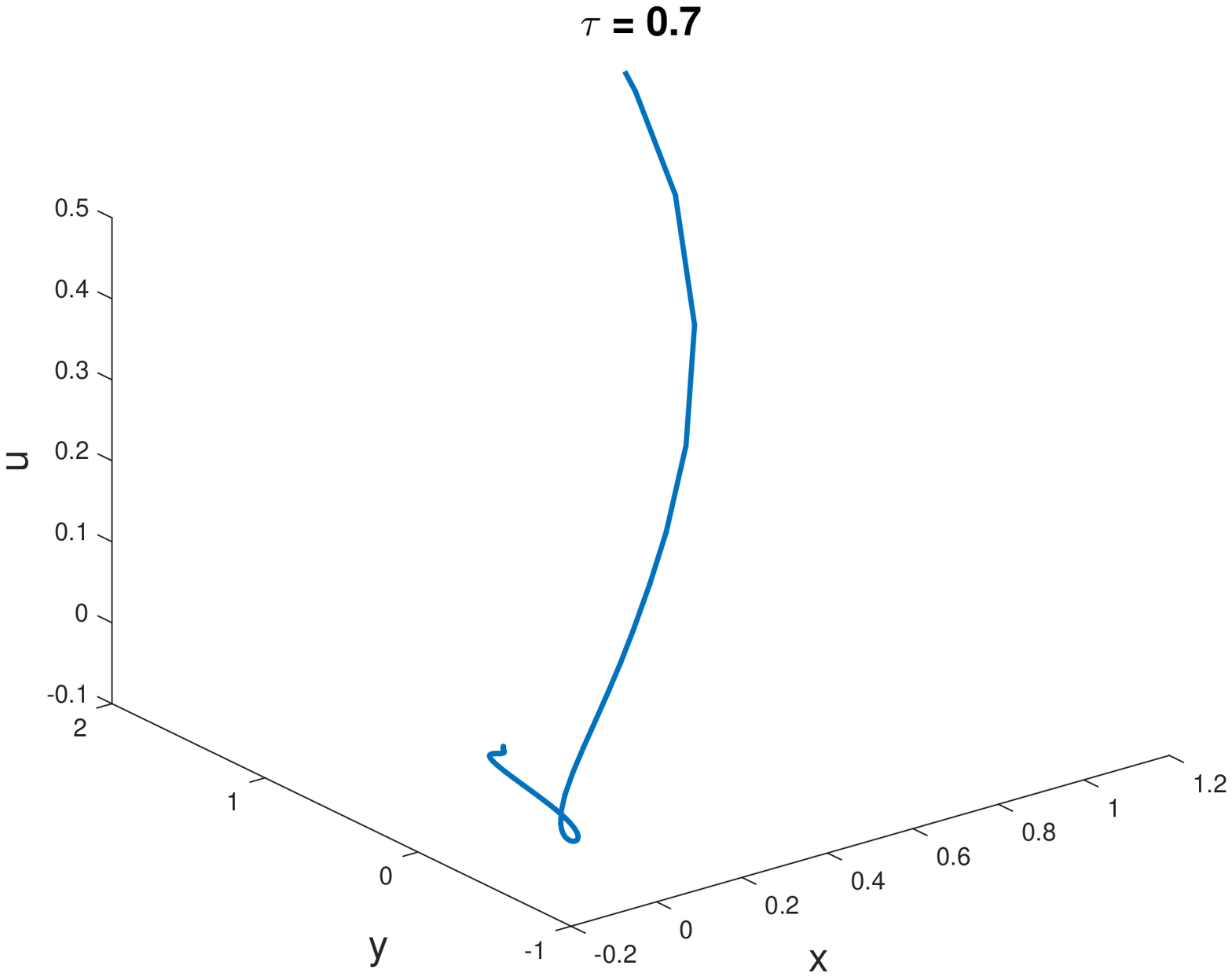}}
\subfigure[]{\label{sek21} 
\includegraphics[scale=0.34]{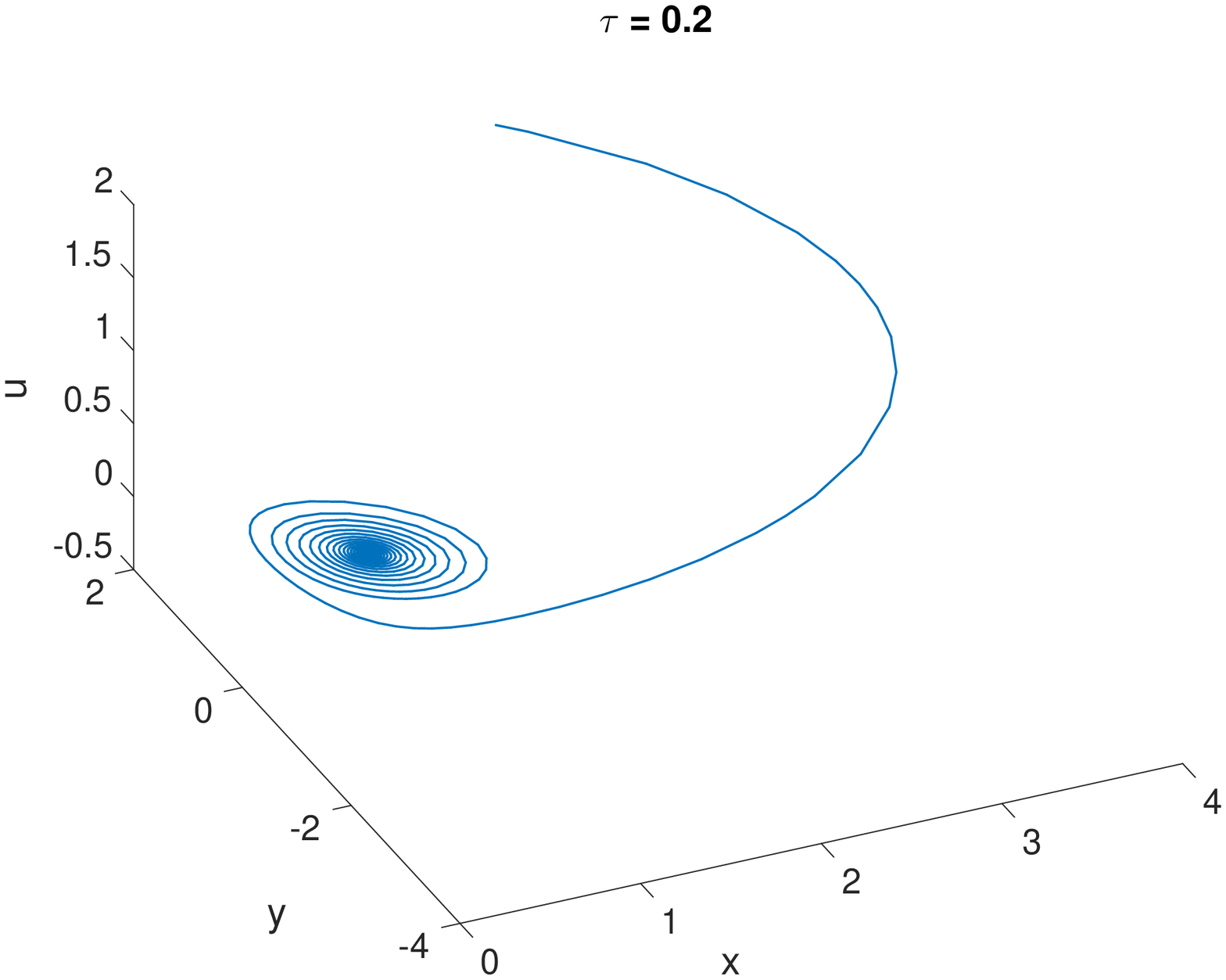}}
\subfigure[]{\label{sek22} 
\includegraphics[scale=0.37]{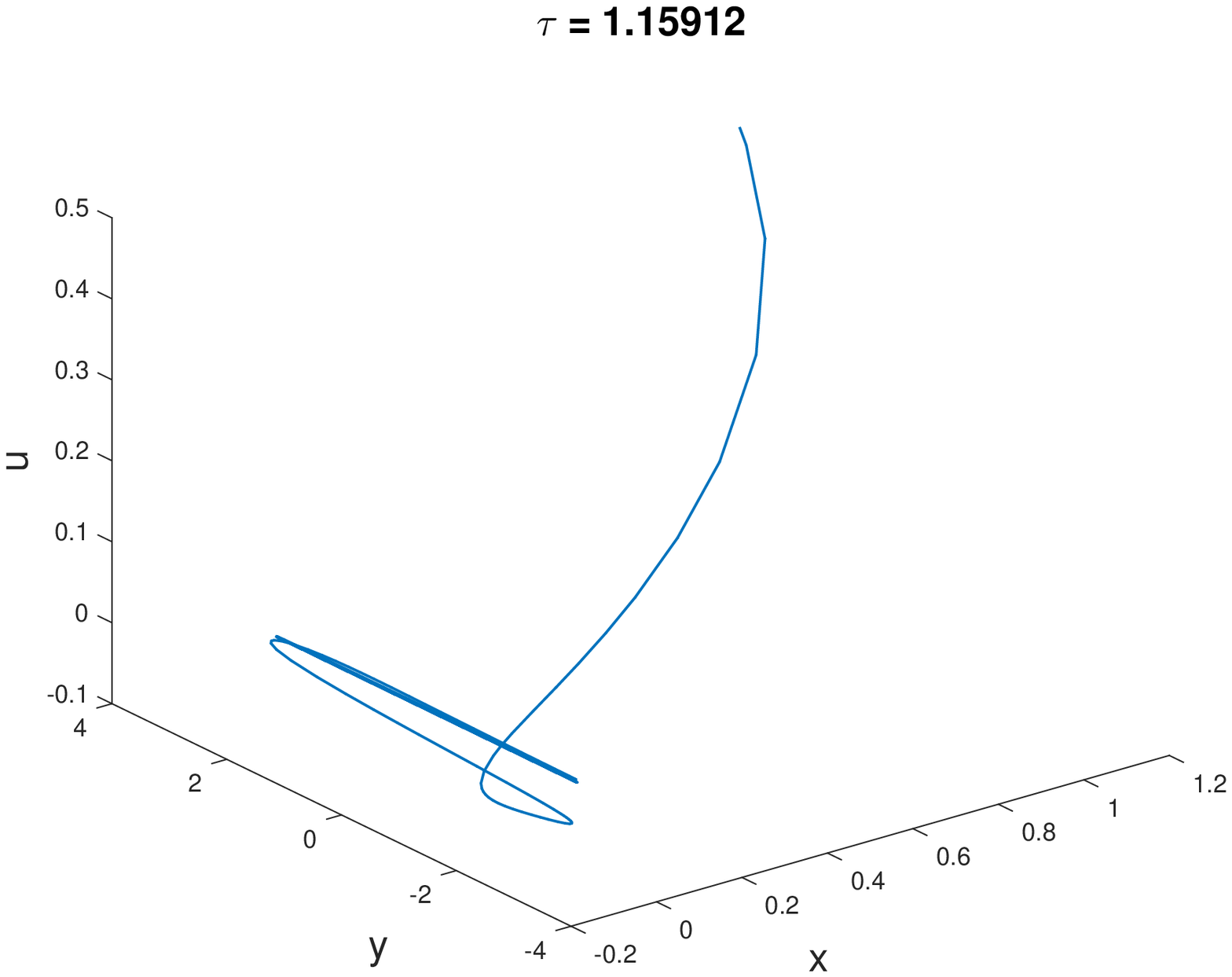}}
\subfigure[]{\label{sek22} 
\includegraphics[scale=0.33]{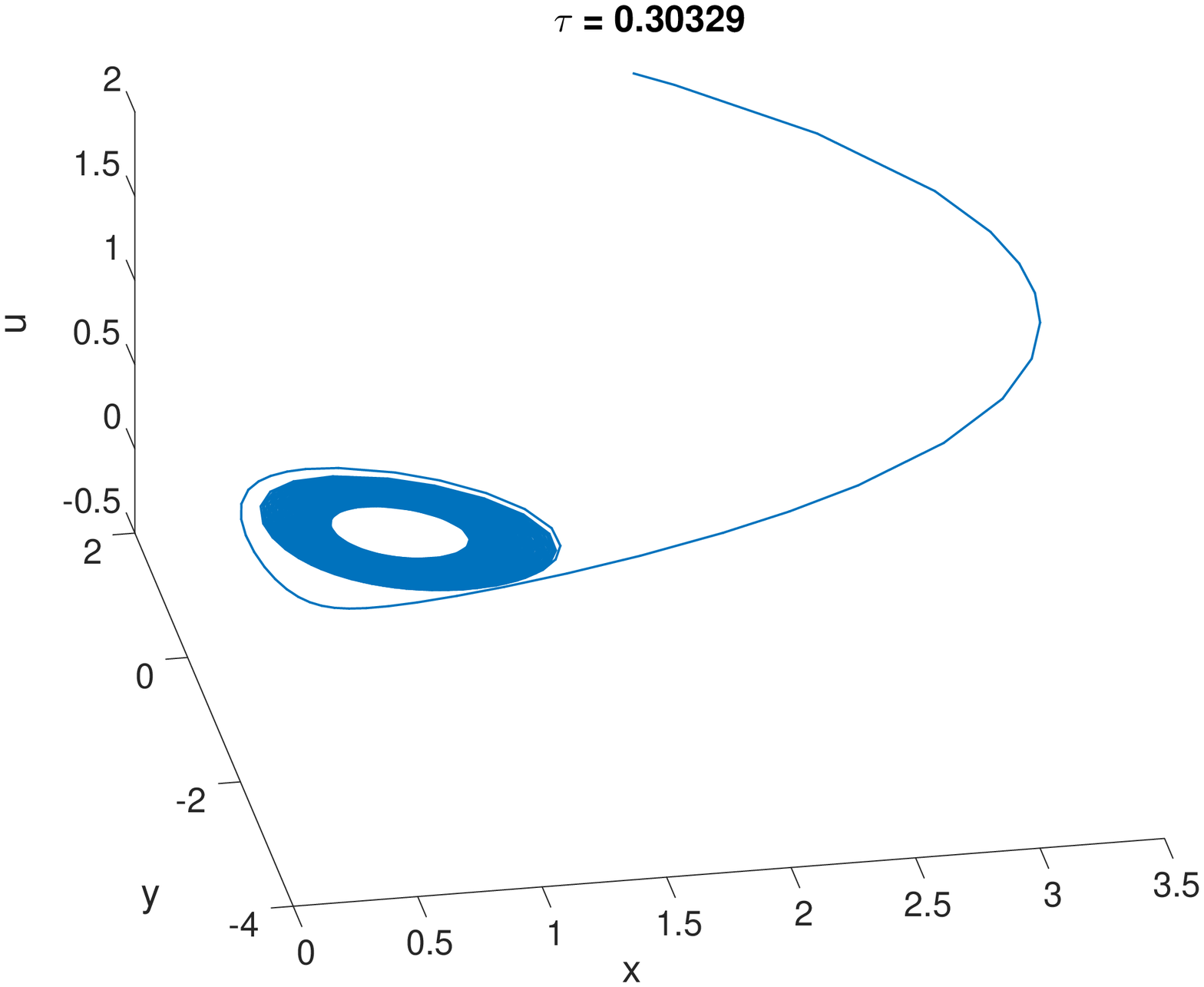}}
\subfigure[]{\label{sek23} 
\includegraphics[scale=0.41]{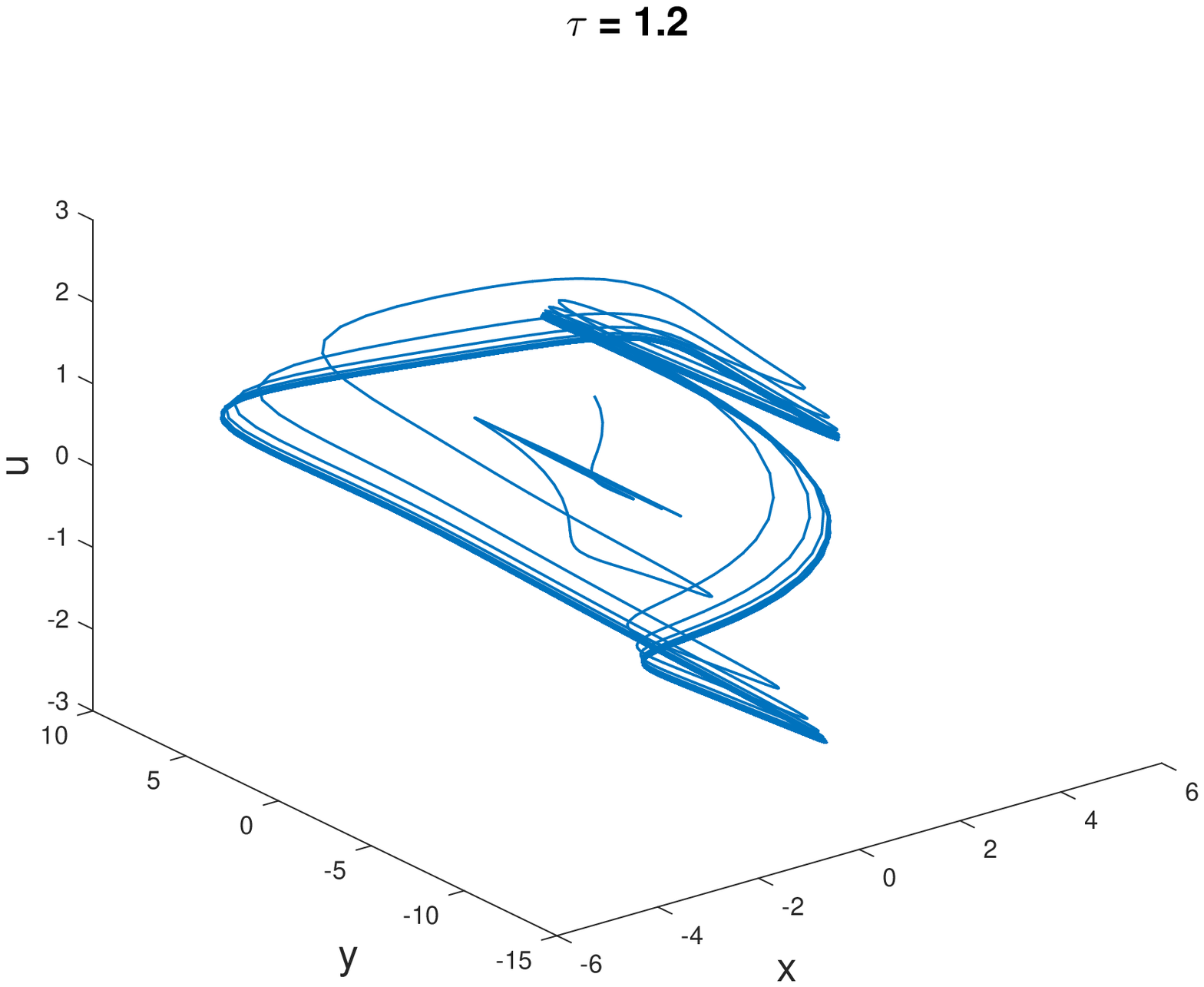}}
\subfigure[]{\label{sek23} 
\hspace{1.5cm}
\includegraphics[scale=0.36]{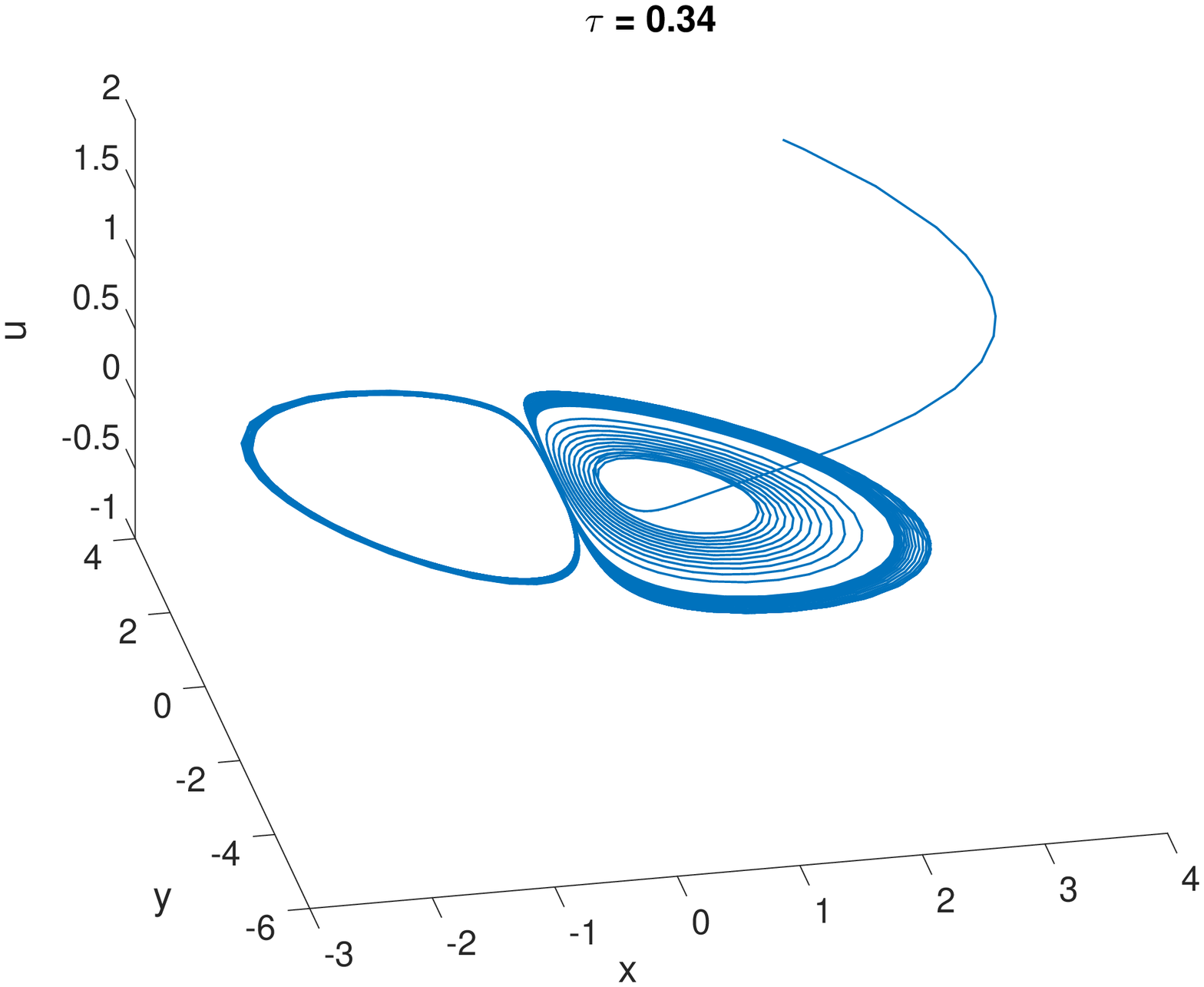}}
\vspace{.1cm}
\caption{Three dimensional phase portraits of the variables $x$, $y$ and $u$ obtained from the numerical solutions of the system \eqref{main} for the chosen parameter values and initial conditions; (a), (c), (e) for $P_0$, (b), (d), (f) for $P_1$.}
\end{figure}


%
%


%
%


\end{document}